# Testing temporal constancy of the spectral structure of a time series

EFSTATHIOS PAPARODITIS

*Department of Mathematics and Statistics, University of Cyprus, P.O. Box 20537, CY-1678 Nicosia, Cyprus. E-mail: stathisp@ucy.ac.cy*

Statistical inference for stochastic processes with time-varying spectral characteristics has received considerable attention in recent decades. We develop a nonparametric test for stationarity against the alternative of a smoothly time-varying spectral structure. The test is based on a comparison between the sample spectral density calculated locally on a moving window of data and a global spectral density estimator based on the whole stretch of observations. Asymptotic properties of the nonparametric estimators involved and of the test statistic under the null hypothesis of stationarity are derived. Power properties under the alternative of a time-varying spectral structure are discussed and the behavior of the test for fixed alternatives belonging to the locally stationary processes class is investigated.

## 1. Introduction

Second order stationarity is an important assumption in the statistical analysis of stochastic processes, allowing for the development of an asymptotic theory capable of investigating properties of many classical statistical inference procedures. However, the assumption that the second order characteristics of a process remain constant over time is often not justified in practice. Many time series show a non-stationary behavior due to smooth changes of their second order structure over time. Based on this observation, interest has been directed toward the analysis of stochastic processes whose spectral characteristics change slowly over time. Priestley (1965) introduced processes with evolutionary spectra and a time-varying spectral representation similar to that of a stationary sequence; cf. also Granger (1964) and Priestley (1988) for an overview. During the last decade, statistical inference for processes showing such a non-stationary behavior has attracted considerable attention in the literature. In this context, asymptotic statistical inference has been made possible by considering arrays of double-indexed processes having time-varying spectral characteristics together with a time rescaling approach which allows for increasing information on the local structure of the underlying non-stationary process as the sample size increases; see Dahlhaus ([1997](#)) and literature on the related concept of locally stationary processes. Mallat *et al.* ([1998](#)) considered adaptive covariance estimators of locally stationary processes, while Nason *et al.* ([2000](#)) considered processes where the time-varying spectral representation and the associated Fourier basis have been replaced by a representation with respect to a wavelet basis. The aforementioned developments in the statistical analysis of locally stationary processes make it important to have powerful statistical tools that detect possible changes in the spectral structure of a time series and evaluate the appropriateness of the preassigned stationary process class.

This paper deals with nonparametric tests of stationarity against the alternative of a time-varying spectral structure. The idea underlying our approach is to compare local sample spectral density estimates based on a moving window of data with a global estimate using the whole set of observations, and to evaluate the difference between the two quantities over the different frequencies and time segments considered using an appropriate $L_2$-type distance measure. To be more specific, the proposed test is based on a smoothed version of the local sample spectral density (local periodogram) rescaled by a global spectral density estimator obtained using the whole stretch of data. For a given time point, the corresponding smoothed statistic squared and integrated over all frequencies is a measure of deviation between the local spectral structure of the observed time series and a global spectral structure associated with the hypothesized stationary behavior. Calculating this local quadratic deviation measure for different instants of time and building a time-averaged version lead to a global measure of stationarity which is used for testing purposes.

To allow for a rigorous asymptotic investigation of the testing methodology proposed under the null and under different non-stationary alternatives, the theoretical framework of locally stationary processes is adopted; cf. Dahlhaus ([1997](#)). The asymptotic behavior of the test statistic under the null of stationarity is derived and its weak convergence to a Gaussian distribution is established. Consistency of the test is proved and its asymptotic power behavior is investigated for fixed alternatives belonging to the class of locally stationary processes. Although some of the results obtained under the null of stationarity can be derived under different assumptions on the underlying stochastic process, the aforementioned theoretical context of locally stationary processes is preferred because it delivers a unified set-up for asymptotic investigations of the test proposed both under the null and under fixed alternatives.

Procedures for testing the constancy of particular characteristics of a process over time have been considered in the literature under various settings and approaches with their main focus directed toward testing the constancy of specific parameters or characteristics in a more or less restrictive parametric set-up. The more general nonparametric problem of testing the constancy of the entire spectral structure of a process without imposing any parametric assumptions on the underlying process class has attracted less interest in the literature and only a few studies of this exist.

Priestley and Rao ([1969](#)) proposed testing the homogeneity of a set of evolutionary spectra evaluated at different instants of time using logarithmic transformations and an analysis of variance framework. For Gaussian processes and for the particular purpose of a change-point detection, Picard ([1985](#)) developed a test based on the difference between (possible time-varying) spectral distribution functions estimated on different parts of the series and evaluated using a supremum-type statistic. Giraitis and Leipus (1992) generalized this approach to the case of linear processes. von Sachs and Neumann (2000) developed a test of stationarity based on empirical wavelet coefficients estimated using localized versions of the periodogram. In a more restrictive context, Sakiyama and Taniguchi (2003) considered testing stationarity versus local stationarity in a parametric set-up. Fitting piecewise autoregressive processes has been also used to detect change-points in the context of time series; Kitagawa and Akaike (1978) used this approach together with AIC-based order selection, while Davis *et al.* (2005) suggested fitting piecewise autoregressive models using a minimum description length procedure together with a generic algorithm to solve the difficult optimization problem.

Compared to the last two approaches mentioned, the method proposed in this paper is entirely nonparametric in that it does not rely on distributional or parametric assumptions on the underlying process class and it is not only restricted on the detection of change-points but also on smooth changes of the spectral structure. Unlike Picard's (1985) proposal focused on change point detection and which is based on the integrated periodogram and supremum-type statistics, our method is based on smooth (kernel) estimates of the (possibly) time-varying spectral density function of the underlying process and uses $L_2$-type distance measures. Compared to the rather heuristic derivations of Priestley and Rao (1969), our set-up allows for a thorough theoretical investigation of the asymptotic properties of the proposed test both under the null and under fixed alternatives. In addition to the aforementioned power investigations, our testing approach differs from that of von Sachs and Neumann (2000) in that it uses global measures of deviations from stationarity, thereby avoiding a multiple testing procedure based on a selected set of estimated wavelet coefficients with conservative critical values obtained via Bonferonni's inequality.

The paper is organized as follows. Section 2 contains the main assumptions imposed on the underlying stochastic process class and introduces the basic statistics used in the sequel. In Section 3, the statistic proposed to test the null hypothesis is presented and its asymptotic behavior under the null of stationarity is investigated. Section 4 deals with the power properties of the test and investigates its behavior for fixed alternatives. Section 5 contains some practical guidelines on how to choose the testing parameters and demonstrates the capability of our method to detect departures from stationarity by means of a numerical example. Proofs of the main results are deferred to Section 6.

## 2. Basic statistics

### 2.1. Preliminary considerations

To provide a theoretical set-up allowing for stochastic processes showing time-varying spectral characteristics and to state rigorously the null and alternative hypothesis of interest, we adopt the framework of locally stationary process (see Dahlhaus (1997) and Dahlhaus and Polonik (2003)) which is flexible enough and covers several interesting cases. Since we restrict ourselves to questions related to the temporal behavior of the spectral structure of a process, it seems natural to consider the general class of time-varying linear processes as the appropriate underlying class of processes.

Consider a triangular array $\{\mathbf{X}_n, n \in \mathbb{N}\}$ of stochastic processes, where $\mathbf{X}_n \equiv \{X_{t,n}, t = 1, 2, \ldots, n\}$ satisfies

***Assumption 1.***

(i) $X_{t,n}$ *has the representation*

$$X_{t,n} = \sum_{j=-\infty}^{\infty} a_{t,n}(j)\varepsilon_{t-j}, \tag{2.1}$$

*where the* $\{\varepsilon_t\}$ *are i.i.d. with* $E(\varepsilon_t) = 0$, $E(\varepsilon_t^2) = 1$, $E(\varepsilon_t^8) < \infty$ *and* $\kappa_4 = \eta - 3$, *where* $\eta = E(\varepsilon_t^4)$;

(ii) $\sup_t |a_{t,n}(j)| \leq K l^{-1}(j)$, *where* $K$ *is a non-negative constant independent of* $n$ *and the positive sequence* $\{l(j), j \in \mathbb{Z}\}$ *satisfies* $\sum_{j=-\infty}^{\infty} |j| l^{-1}(j) < \infty$;

(iii) *functions* $a(\cdot, j) : (0, 1] \to \mathbb{R}$ *with*

$$\sup_{u \in [0,1]} \left| \frac{\partial^2 a(u, j)}{\partial u^2} \right| \leq \frac{K}{l(j)} \tag{2.2}$$

*exist such that*

$$\sup_{1 \leq t \leq n, n \in \mathbb{N}} \left| a_{t,n}(j) - a\left(\frac{t}{n}, j\right) \right| \leq \frac{K}{n l(j)}. \tag{2.3}$$

Let $A_{t,n}(\lambda) = \sum_{j=-\infty}^{\infty} a_{t,n}(j) \exp\{-\mathrm{i}\lambda t\}$ and denote by $f(u, \lambda)$ the time-varying spectral density $f(u, \lambda) = (2\pi)^{-1} |A(u, \lambda)|^2$, where $A(u, \lambda) = \sum_{j=-\infty}^{\infty} a(u, j) \exp\{-\mathrm{i}j\lambda\}$. Note that by the smoothness properties of $a(\cdot, j)$, the triangular array $\{\mathbf{X}_n, n \in \mathbb{N}\}$ uniquely determines the local spectral density $f(u, \lambda)$, and that (2.3) implies

$$\sup_{t, \lambda} \left| \frac{1}{2\pi} A_{t,n}(\lambda) \overline{A_{t,n}(\lambda)} - f(t/n, \lambda) \right| = \mathrm{O}(n^{-1}). \tag{2.4}$$

Let $\mathcal{F}_{LS}$ be the class of stochastic processes satisfying Assumption 1 and let $\mathcal{F}_S$, $\mathcal{F}_S \subset \mathcal{F}_{LS}$, be the subclass of stationary processes, that is,

$$\mathcal{F}_S = \left\{ \mathbf{X}_n : X_{t,n} = \sum_{j=-\infty}^{\infty} a(j)\varepsilon_{t-j} \ \forall t \text{ and } n \in \mathbb{N} \text{ and } \sum_{j=-\infty}^{\infty} |j||a(j)| < \infty \right\}.$$

Given observations $X_{1,n}, X_{2,n}, \ldots, X_{n,n}$, the problem considered is that of testing the hypothesis that the spectral structure of the underlying process is constant over time. To precisely define the corresponding null and alternative hypotheses, consider for $\lambda \in [-\pi, \pi]$ the time-averaged local spectral density $g(\lambda)$ defined by

$$g(\lambda) = \int_0^1 f(u, \lambda) \, \mathrm{d}u.$$

By the smoothness properties of $f(\cdot, \lambda)$, $g(\lambda)$ is well defined and has all properties in common with the spectral density of a real-valued stationary process, that is, it is symmetric, non-negative definite and satisfies $\int_{-\pi}^{\pi} g(\lambda)\,\mathrm{d}\lambda < \infty$. Furthermore, $f(u, \cdot) = g(\cdot)$ if and only if for every $\lambda \in [-\pi, \pi]$, $f(u, \lambda)$ is an a.e. constant function of the time variable $u \in [0, 1]$. We require that $g$ also satisfies the following assumption.

**Assumption 2.**

$$\inf_{\lambda \in [-\pi, \pi]} g(\lambda) > 0. \tag{2.5}$$

Based on the previous considerations, the testing problem investigated is described by the following pair of null and alternative hypotheses

$$H_0 : f(u, \cdot) = g(\cdot) \qquad \text{a.e.}$$

$$\text{vs.} \tag{2.6}$$

$$H_1 : \lambda(\{u : f(u, \cdot) \neq g(\cdot), u \in [0, 1]\}) > 0,$$

where $\lambda(A)$ denotes the Lebesgue measure of the set $A \in \mathcal{B}([0, 1])$.

Let $I_n(\lambda)$, $\lambda \in [-\pi, \pi]$ be the periodogram calculated using the whole time series, that is,

$$I_n(\lambda) = \frac{1}{2\pi n} \left| \sum_{t=1}^{n} X_{t,n} \mathrm{e}^{-\mathrm{i}\lambda t} \right|^2,$$

and let $K$ be a kernel function $K : \mathbb{R} \to \mathbb{R}$ satisfying the following assumption.

**Assumption 3.** $K$ *is Lipschitz continuous, bounded and symmetric with support* $[-\pi, \pi]$ *such that* $(2\pi)^{-1} \int_{-\pi}^{\pi} K(x)\,\mathrm{d}x = 1$.

Consider the smoothed periodogram

$$\widehat{g}_h(\lambda) = \frac{1}{n} \sum_{j \in \mathbb{Z}} K_h(\lambda - \omega_j) I_n(\omega_j), \tag{2.7}$$

where $\omega_j = 2\pi j/n$, $j = 0, \pm 1, \pm 2, \ldots$, $K_h(\cdot) = h^{-1} K(\cdot/h)$ is the scaled kernel and $h > 0$ the bandwidth. It is well known that if the null hypothesis $H_0$ is true, then

$$\sup_{\lambda \in [-\pi, \pi]} |\widehat{g}_h(\lambda) - f(\lambda)| \to 0, \tag{2.8}$$

in probability, where $f(\lambda) = (2\pi)^{-1} \sigma^2 |A(\mathrm{e}^{-\mathrm{i}\lambda})|^2$, $A(\mathrm{e}^{-\mathrm{i}\lambda}) = \sum_{j=-\infty}^{\infty} a(j) \exp\{-\mathrm{i}j\lambda\}$, denotes the time-invariant spectral density of the underlying stationary process. The above result can also be extended to the case where the alternative $H_1$ is true, that is, where $\mathbf{X}_n$ is locally stationary with a time-varying spectral structure in the sense of Assumption 1. In particular, as Lemma 6.1 shows, we have, in this case,

$$\sup_{\lambda \in [-\pi, \pi]} \left| \widehat{g}_h(\lambda) - \int_0^1 f(u, \lambda) \, \mathrm{d}u \right| \to 0 \tag{2.9}$$

in probability, that is, for locally stationary processes, $\widehat{g}_h(\lambda)$ is a uniformly consistent estimator of the time-averaged local spectral density $g(\lambda)$.

## 2.2. Smoothed rescaled local periodogram

To construct out test statistic, let $m_n$ be a positive integer such that $0 < m_n < n$ and consider for $\lambda \in [-\pi, \pi]$ and $u \in (0, 1)$ the local periodogram $I_{m_n}(u, \lambda)$ calculated using a subset of $m_n$ tapered observations around location $u$, that is,

$$I_{m_n}(u, \lambda) = \frac{1}{2\pi H_{2, m_n}(0)} \left| \sum_{t=1}^{m_n} h_{t, m_n} X_{t + [un] - M_n - 1, n} \mathrm{e}^{-\mathrm{i}\lambda t} \right|^2,$$

where $M_n = m_n/2$ and $u \in (0, 1)$ such that $M_n + 1 \leq [un] \leq n - M_n + 1$.

Furthermore, $h_{t,n} = h(t/n)$ is a taper function and

$$H_{k,n}(\lambda) = \sum_{s=1}^{n} h(s/n)^k \mathrm{e}^{-\mathrm{i}\lambda s}$$

its Fourier transform. Write $H_n(\cdot)$ for $H_{1,n}(\cdot)$.

**Assumption 4.** *The taper function* $h : \mathbb{R} \to [0, 1]$ *is of bounded variation and vanishes outside the interval* $[0, 1]$.

Note that the use of taper reduces the bias of the local periodogram due to the well-known leakage effect and it simultaneously reduces the bias due to the (possible) non-stationarity of the series. For the asymptotic theory developed in this paper, the use of such a taper is necessary in order to control the bias of the local periodogram, under both the null and the alternative. Note that we do not use a data taper to calculate the global periodogram $I_n(\lambda)$ based on the whole stretch of observations mainly because we would then get

$$E(I_n(\lambda_k)) = H_{2,n}(0)^{-1} \sum_{t=1}^{n} h^2(t/n) f(t/n, \lambda_k) + \mathrm{O}(n^{-1} \log(n)),$$

which implies that under the assumptions made, the smoothed global periodogram $\widehat{g}_h(\lambda)$ does not converge under the alternative to the desired limit $g(\lambda) = \int_0^1 f(u, \lambda) \, \mathrm{d}u$.

Now, the basic statistic used in the sequel is the following kernel-smoothed version of the rescaled local periodogram

$$V_n(u,\lambda) = \frac{1}{m_n} \sum_{j\in\mathbb{Z}} K_b(\lambda - \lambda_j)\left(\frac{I_{m_n}(u,\lambda_j)}{\widehat{g}_h(\lambda_j)} - 1\right), \qquad (2.10)$$

where $\lambda_j = 2\pi j/m_n$, $j = 0, \pm 1, \pm 2, \ldots$. To understand the motivation leading to $V_n(u,\lambda)$, consider

$$\widetilde{V}_n(u,\lambda) = \frac{1}{m_n} \sum_{j\in\mathbb{Z}} K_b(\lambda - \lambda_j)\left(\frac{I_{m_n}(u,\lambda_j)}{g(\lambda_j)} - 1\right), \qquad (2.11)$$

which differs from $V_n(u,\lambda)$ by the fact that the true function $g(\lambda)$ is used. Under the assumptions made, and if the null hypothesis is true, then $\widetilde{V}_n(u,\lambda)$ is a nonparametric (kernel) estimator of the mean function $E(I_{m_n}(u,\lambda)/f(\lambda) - 1) = \mathrm{O}(m_n^{-2})$ that converges to zero uniformly in $u$ and $\lambda$ as $m_n \to \infty$. Thus, under the null hypothesis, we have, for all $u \in (0,1)$ and $\lambda \in [-\pi, \pi]$, that

$$\widetilde{V}_n(u,\lambda) \to 0$$

in probability. On the other hand, if the alternative is true, that is, if the spectral structure of the underlying process varies over time, then it follows using Lemma 6.1, that

$$\widetilde{V}_n(u,\lambda) \to \left(\frac{f(u,\lambda)}{\int_0^1 f(u,\lambda)\,\mathrm{d}u} - 1\right)$$

in probability, where the limiting function on the right-hand side is different from the zero function because $f(u,\cdot) = \int_0^1 f(u,\cdot)\,\mathrm{d}u$ if and only if $f(u,\cdot)$ is an a.e. constant function of the time parameter $u \in [0,1]$. Thus, under the null hypothesis of stationarity, we expect $\widetilde{V}_n(\cdot,\cdot)$ to be close to the zero function, while this will not be the case under the alternative.

Note that for a given $u \in (0,1)$, the statistic (2.10) has the same asymptotic behavior as the statistic $\widetilde{V}_n(\cdot)$ in which the true function $g$ appears. This is true because (2.10) is based on observations within a window of length $m_n$ while $\widehat{g}_h$ is calculated using the whole stretch of $n$ data points. In fact, $V_n(u,\lambda) = \widetilde{V}_n(u,\lambda) + \mathrm{O}_P(\sup_{\lambda_j} |\widehat{g}_h(\lambda_j) - g(\lambda_j)|)$. For statistics based on time-averaged versions of (2.10), such as those considered in the sequel, we show that under appropriate conditions on the behavior of the smoothing parameters involved, the effects of using the nonparametric estimator $\widehat{g}_h$ are asymptotically negligible.

# 3. Testing the null hypothesis

## 3.1. The test statistic

To make the aforementioned behavior of the statistic $V_n(u, \lambda)$ useful for testing the null hypothesis of stationarity, let $0 < u_1 < u_2 < \cdots < u_N < 1$ be a set of $N = N(n) \in \mathbb{N}$ distinct points in the interval $(0, 1)$. Set $u_0 \equiv 0$ and let $d_{j,n} = u_j - u_{j-1}$, $j = 1, 2, \ldots, N$ be the distance between two consecutive points. Taking equidistant points, that is, $d_{j,n} = d_n = (N+1)^{-1}$, a choice of $u_j$ could be $u_j = t_j/n$, where $t_j = S(j-1) + m_n/2$ for $j = 1, 2, \ldots, N$. The positive integer $S = S(n)$ denotes the shift from segment to segment while $n = S(N-1) + m_n$. Using a squared deviation criterion, a useful approach to test the null hypothesis is to average over the different time points and over the different frequencies the squared statistic $V_n(t_s/n, \lambda)$. This leads to the test statistic

$$T_n = \frac{1}{N} \sum_{s=1}^{N} \int_{-\pi}^{\pi} V_n^2\left(\frac{t_s}{n}, \lambda\right) \mathrm{d}\lambda \tag{3.1}$$

considered in the sequel. Note that due to averaging over $N$ distinct time points, $T_n$ is an overall measure of deviation of the local spectral structure of the process from its hypothesized global behavior under the null hypothesis.

To investigate the asymptotic behavior of $T_n$ under the null of stationarity, the following assumptions are imposed.

**Assumption 5.** *As $n \to \infty$,*

  (i) $h \sim n^{-\lambda}$ *for some* $3/20 < \lambda < 1/3$;

  (ii) $b \sim m_n^{-\lambda}$;

  (iii) $S_n = [m_n/c]$, *where $c$ is a fixed positive integer such that $1 \le c < m_n$;*

  (iv) $m_n \sim n^\delta$ *for some* $\delta \in (\delta_1, \delta_2)$, *where* $\delta_1 = \max\{1/(3-\lambda), \lambda/(1-\lambda)\}$ *and* $\delta_2 = \min\{(8\lambda - 1)/(1-\lambda), (1-2\lambda)/(1-\lambda)\}$.

Some comments concerning the above assumptions are in order. Part (i) and (ii) ensure consistency of the corresponding global and local spectral estimators. For simplicity of derivations, we allow both smoothing bandwidths $h$ and $b$ to converge to zero at the same rate, although different rates are also possible. Part (iii) puts some conditions on the behavior of the shift $S_n$ between the different segments used to calculate the statistic $V_n^2(u_s, \lambda)$. The condition $m_n \ge S_n$ seems reasonable since it makes no sense to omit data ($m_n < S_n$). Although the requirement for a fixed $1 \le c < m_n$ is not restrictive in practice, it is necessary from a technical point of view in order to control the degree of overlapping between the segments used. Part (iv) implies that $m_n/n \to 0$ as $n \to \infty$, from which it follows, because of $m_n = cS_n$ and $NS_n = \mathrm{O}(n)$, that the number $N$ of time points $u_j = t_j/n$ at which $V_n(u_j, \lambda)$ is calculated increases to infinity as $n \to \infty$. This is important for a good power behavior of the test since, in this case, the (rescaled) distance $d_n = S_n/n$ between the time points $t_j/n$ in the interval $(0, 1)$ over which the behavior of $V_n(t_j/n, \lambda)$ is evaluated goes to zero as $n \to \infty$; cf. Section 4. However, part (iv) introduces some restrictions regarding the rates at which the time window $m_n$ and, consequently, the number of time points $N$, increase to infinity as the sample size $n$ increases. Essentially, these conditions prevent both $m_n$ and $N$ from increasing too fast with respect to $n$ and seem necessary in order to make the effects of estimating the time-integrated spectral density $g$ on the distribution of the test statistic $T_n$ asymptotically negligible.

**Theorem 3.1.** *Let $H_0$ be true and Assumptions 1–5 be satisfied. As $n \to \infty$,*

$$m_n\sqrt{Nb}\,T_n - \mu_n \Rightarrow \mathcal{N}(0, \tau^2),$$

*where*

$$\mu_n = M_H\left[\sqrt{\frac{N}{b}}\int_{-\pi}^{\pi} K^2(x)\,\mathrm{d}x + \sqrt{Nb}\left(\frac{1}{4\pi}\int_{-2\pi}^{2\pi}(K*K)(y)\,\mathrm{d}y + 2\pi\kappa_4\right)\right],$$

$$\tau^2 = \frac{2}{\pi^2}M_{c,H}^2\int_{-2\pi}^{2\pi}(K*K)^2(y)\,\mathrm{d}y,$$

$$M_{c,H}^2 = H_2^{-4}\sum_{|s|<c}\left(\int_{[0,1-|s|/c]}h^2(x)h^2(x+|s|/c)\,\mathrm{d}x\right)^2,$$

$M_H = H_4/H_2^2$, $H_k = \int_0^1 h^k(u)\,\mathrm{d}u$ *and $K*K(\cdot)$ denotes convolution of the kernel $K$.*

The fact that the centering sequence $\mu_n$ in the above theorem depends on the fourth order cumulant $\kappa_4$ of the error process is due to time-averaging the weak and $\mathrm{O}(m_n^{-1})$-vanishing covariance of the local periodogram ordinates at different Fourier frequencies. To implement the statistic $T_n$ in practice, a consistent estimator of $\kappa_4$ is needed. Such a nonparametric estimator is given by

$$\widehat{\kappa}_4 = \left(2\pi\widehat{g}_{2,h}(0) - 4\pi\int_{-\pi}^{\pi}\widehat{g}_h^2(\lambda)\,\mathrm{d}\lambda\right)\bigg/\left(\int_{-\pi}^{\pi}\widehat{g}_h(\lambda)\,\mathrm{d}\lambda\right)^2, \tag{3.2}$$

where $\widehat{g}_{2,h}(\cdot)$ is the estimated spectral density of the squared process $X_{t,n}^2$; cf. Janas and Dahlhaus (1994) and Grenander and Rosenblatt (1956). Based on Theorem 3.1, an asymptotically $\alpha$-level test, $\alpha \in (0,1)$, is obtained by rejecting the null hypothesis if $m_n\sqrt{Nb}\,T_n \geq \widehat{\mu}_n + \tau z_\alpha$, where $z_\alpha$ denotes the upper $\alpha$-percentile of the standard Gaussian distribution and $\widehat{\mu}_n$ is the estimator of $\mu_n$ obtained by replacing $\kappa_4$ by (3.2).

# 4. Power properties

We begin our power investigations with the following theorem which establishes consistency of the proposed test.

**Theorem 4.1.** *Suppose that the triangular array* $\{\mathbf{X}_n, n \in \mathbb{N}\}$ *possesses a local spectral density* $f(u, \lambda)$, $f(\cdot, \cdot) \in L_2([0, 1] \times [-\pi, \pi])$. *Assume that the set* $\{u : u \in [0, 1], f(u, \lambda) \neq g(\lambda)\} \subset [0, 1]$ *has positive Lebesgue measure and that the following conditions are satisfied:*

$$\widehat{g}_h(\lambda) \underset{n \to \infty}{\to} g(\lambda) \qquad \text{in probability}, \qquad \limsup_{n \to \infty} \int_{-\pi}^{\pi} \widehat{g}_h(\lambda) \, \mathrm{d}\lambda \leq \int_{-\pi}^{\pi} g(\lambda) \, \mathrm{d}\lambda$$

*and*

$$\int_{-\pi}^{\pi} \left| m_n^{-1} \sum_{j \in \mathbb{Z}} K_b(\lambda - \lambda_j) I_{m_n}(u, \lambda_j) - f(u, \lambda) \right| \mathrm{d}\lambda = \mathrm{O}_P(r_n)$$

*uniformly in* $u$, *where* $r_n$ *is a zero sequence. Then, as* $n \to \infty$,

$$T_n \to D^2 = \int_0^1 \int_{-\pi}^{\pi} \left( \frac{f(u, \lambda)}{g(\lambda)} - 1 \right)^2 \mathrm{d}\lambda \, \mathrm{d}u$$

*in probability.*

Notice that the above conditions are satisfied if, for instance, the process under the alternative is a locally stationary process in the sense of Assumption 1. Observe further that the limit $D^2$ is an $L_2$-measure of deviation of the locally stationary process with time-varying local spectral density $f(u, \lambda)$ from the null, that is, from a stationary process with spectral density $g(\lambda) = \int_0^1 f(u, \lambda) \, \mathrm{d}u$. Recall that $g(\lambda)$ is a proper spectral density and note that $g(\lambda)$ can be interpreted as the time-invariant spectral density which best approximates the time-varying local spectral density $f(u, \lambda)$ in the sense that for each frequency $\lambda \in [-\pi, \pi]$, the function $g(\lambda) = \int_0^1 f(u, \lambda) \, \mathrm{d}u$ satisfies $g(\lambda) = \operatorname{argmin}_c \int_0^1 (f(u, \lambda) - c)^2 \, \mathrm{d}u$. Positive values of $D^2$ are tantamount to a deviation of the locally stationary process $\{\mathbf{X}_n\}$ from its best approximating "stationary counterpart", where the latter can be understood as a zero-mean stationary process with autocovariance function $\gamma(\cdot) : \mathbb{Z} \to \mathbb{R}$ given by $\gamma(h) = \int_{-\pi}^{\pi} g(\lambda) \exp\{\mathrm{i}h\lambda\} \, \mathrm{d}\lambda$.

## 4.1. Locally stationary alternatives

Although Theorem 4.1 provides some useful insights into the consistency properties of the proposed test, an asymptotic analysis of its behavior for fixed alternatives is very informative and important since it leads to useful approximations of its power function and of the probability of the type II error. In the following, we investigate the power behavior of the test for fixed alternatives that belong to the class of locally stationary process. For this, the following assumption is imposed.

**Assumption 6.** *As* $n \to \infty$,

(i) $h \sim n^{-\lambda}$ *for some* $1/4 < \lambda < 1/3$;
(ii) $b \sim m_n^{-\lambda}$;

(iii) $S_n = [m_n/c]$, where $c$ is a fixed integer such that $1 \le c < m_n$;

(iv) $m_n \sim n^\delta$ for some $\delta$, $1/2 < \delta < \min\{(8\lambda - 1)/(1 - \lambda), (1 - 2\lambda)/(1 - \lambda)\}$.

The lower bound in part (i) above implies that the bandwidth used to obtain the nonparametric estimators of the time-integrated function $g(\lambda)$ and of the rescaled local spectral density $f(u, \lambda)/g(\lambda)$ is relative small, that is, the bias of the corresponding estimators is of a smaller order of magnitude than their variance. Relaxing this assumption to $3/20 < \lambda < 1/3$ (cf. Assumption 5) will introduce an additional bias term to the centering sequence of the asymptotic distribution given in Theorem 4.2 below.

**Theorem 4.2.** *Suppose that* $\mathbf{X}_n \in \mathcal{F}_{LS} \setminus \mathcal{F}_S$, *that Assumptions 1–4 and Assumption 6 are satisfied and that* $E|\varepsilon_1|^k < \infty$ *for all* $k \in \mathbb{N}$. *If* $n \to \infty$, *then*

$$\sqrt{m_n N}(T_n - \delta_n^2 - D_n^2) \Rightarrow \mathcal{N}(0, v^2),$$

*where*

$$\delta_n^2 = \frac{1}{2\pi m_n b} M_H \int_{-\pi}^{\pi} K^2(x)\,\mathrm{d}x \int_0^1 \int_{-\pi}^{\pi} \frac{f^2(u, \lambda)}{g^2(\lambda)}\,\mathrm{d}\lambda\,\mathrm{d}u,$$

$$D_n^2 = \frac{1}{N} \sum_{s=1}^{N} \int_{\pi}^{\pi} \left(\frac{f(u_s, \lambda)}{g(\lambda)} - 1\right)^2 \mathrm{d}\lambda$$

*and* $v^2 = 8\pi(v_1 + L_{c,H}^2 v_2 - 2v_3)$, *with*

$$v_1 = 2 \int_{-\pi}^{\pi} v^2(\lambda) g^2(\lambda)\,\mathrm{d}\lambda + \int_{-\pi}^{\pi} \int_{-\pi}^{\pi} v(\lambda_1) v(\lambda_2) g_4(\lambda_1, -\lambda_1, \lambda_2)\,\mathrm{d}\lambda_1\,\mathrm{d}\lambda_2,$$

$$v_2 = 2 \int_0^1 \int_{-\pi}^{\pi} w^2(u, \lambda) f^2(u, \lambda)\,\mathrm{d}\lambda\,\mathrm{d}u + \kappa_4 \int_0^1 \left[\int_{-\pi}^{\pi} w(u, \lambda) f(u, \lambda)\,\mathrm{d}\lambda\right]^2 \mathrm{d}u,$$

$$v_3 = 2 \int_0^1 \int_{-\pi}^{\pi} v(\lambda) w(u, \lambda) f^2(u, \lambda)\,\mathrm{d}\lambda\,\mathrm{d}u$$

$$\quad + \kappa_4 \int_0^1 \int_{-\pi}^{\pi} \int_{-\pi}^{\pi} v(\lambda_1) w(u, \lambda_2) f(u, \lambda_1) f(u, \lambda_2)\,\mathrm{d}\lambda_1\,\mathrm{d}\lambda_2\,\mathrm{d}u,$$

$$w(u, \lambda) = \frac{1}{g(\lambda)}\left(\frac{f(u, \lambda)}{g(\lambda)} - 1\right), \qquad v(\lambda) = \int_0^1 \frac{f(u, \lambda)}{g^2(\lambda)}\left(\frac{f(u, \lambda)}{g(\lambda)} - 1\right)\mathrm{d}u,$$

$M_H = H^4/H_2^2$, $L_{c,H}^2 = H_2^{-2} \sum_{|s| < c} \int_{[0, 1-|s|/c]} h^2(x) h^2(x + |s|/c)\,\mathrm{d}x$, $g_4(\lambda_1, \lambda_2, \lambda_3) = \int_0^1 f_4(u, \lambda_1, \lambda_2, \lambda_3)\,\mathrm{d}u$ and $f_4(u, \lambda_1, \lambda_2, \lambda_3) = (2\pi)^{-3}\kappa_4 A(u, \lambda_1) A(u, \lambda_2) A(u, \lambda_3) A(u, \sum_{s=1}^3 \lambda_s)$ *the fourth order local cumulant spectrum.*

By the above theorem, asymptotic normality of the test statistic $T_n$ still holds under fixed alternatives, but with a different rate of convergence than under the null of stationar-ity. In fact, while under the null hypothesis the variance of $T_n$ is of order $\mathrm{O}(m_n^{-2}N^{-1}b^{-1})$, under fixed alternatives, the variance of the test statistic is of order $\mathrm{O}(m_n^{-1}N^{-1})$.

Theorem 4.2 is useful for calculating, by means of numerical integration, the asymptotic power of the test for fixed alternatives belonging to the locally stationary process class considered. This is important if, for instance, one is interested in quantifying the asymptotic probability of the type II error associated with using a stationary model. In particular, if $f(u, \lambda)$ is the true local spectral density, then the probability of rejecting the null hypothesis is approximately equal to

$$
\begin{aligned}
P(H_0 \text{ is rejected}) &\approx 1 - \Phi\left(\frac{1}{v}\left(\frac{\tau z_\alpha + \mu_n}{\sqrt{m_n b}} - \sqrt{m_n N}D_n^2 - \sqrt{m_n N}\delta_n^2\right)\right) \\
&\approx 1 - \Phi\left(\frac{1}{v}\left(\frac{\tau z_\alpha}{\sqrt{m_n b}} - \sqrt{m_n N}D_n^2 - \frac{M_H}{2\pi}\frac{\sqrt{N}}{\sqrt{m_n b}}D^2\right)\right) \quad (4.1) \\
&\approx 1 - \Phi\left(-\frac{\sqrt{m_n N}}{v}D_n^2 - \frac{M_H}{2\pi v}\frac{\sqrt{N}}{\sqrt{m_n b}}D^2\right),
\end{aligned}
$$

where $v = \sqrt{v^2}$, with $v^2$, $\delta_n^2$ and $D_n^2$ as given in Theorem 4.2 and $\mu_n$ and $\tau$ as given in Theorem 3.1. The second approximation above follows using $\mu_n/\sqrt{m_n b} - \sqrt{m_n N}\delta_n^2 = \sqrt{N}/(\sqrt{m_n b})M_H(2\pi)^{-1}D^2 + \mathrm{O}(N^{1/2}m_n^{-1/2})$. Notice that the power of the test is dominated by the term $D_n^2$ which is a discrete-time version of the quadratic deviation $D^2$ given in Theorem 4.1 and that $\sqrt{m_n N}D_n^2$ goes to infinity at a $\sqrt{n}$-rate as $n \to \infty$.

# 5. Applications

## 5.1. Some guidelines for choosing the testing parameters

While the results of the proposed testing procedure seem to be less sensitive with respect to the particular choice of the taper function $h(\cdot)$ and the kernel $K(\cdot)$, the choice of the other tuning parameters is more influential and more difficult to solve from a theoretical point of view. Based on some rather heuristic considerations, in this section, we provide some guidelines on how to choose these parameters in practice.

Concerning the smoothing bandwidths $h$ and $b$, one approach is to choose these parameters in a way which leads to "good estimates" of the unknown functions $g(\lambda)$ and $f(u, \lambda)/g(\lambda)$, respectively. A procedure for selecting the smoothing bandwidth in the context of spectral density estimation for stationary process is offered by the cross-validation criterion proposed by Beltrão and Bloomfield ([1987](#)); see also Robinson ([1991](#)). Adapted to our context and under the null hypothesis, such a rule implies that $h$ can be chosen as the minimizer of the objective function

$$
CV(h) = \frac{1}{L_n}\sum_{j=1}^{L_n}\left\{\log(\widehat{g}^{-j}(\omega_j)) + \frac{I_n(\omega_j)}{\widehat{g}^{-j}(\omega_j)}\right\},
$$

where $\widehat{g}^{-j}(\omega_j)$ is the leave-out-$j$ version of the estimator $\widehat{g}_h(\omega_j)$ given by $\widehat{g}^{-j}(\omega_j) = n^{-1} \times \sum_{s \in L_{n,j}} K_h(\omega_j - \omega_s) I_n(\omega_s)$ and $L_{n,j} = \{s : -L_n \leq s \leq L_n \text{ and } j - s \neq \pm j \bmod L_n\}$. Although the above approach modified appropriately can also be applied to select the bandwidth $b$ used within each segment, a computationally simpler rule can be obtained via the following considerations. Suppose that the global bandwidth $h$ behaves under the null of stationarity like $h = C n^{-\lambda}$, where $\lambda$ is known and the constant $C$ depends on characteristics of the unknown function $g(\lambda) = f(u, \lambda)$. Suppose, further, that the bandwidth $b$ has a similar behavior to $h$, but adjusted for the smaller number of observations, that is, $b = C m_n^{-\lambda}$. Then, for $h$ chosen, $b$ can be selected according to the rule $b = h(n/m_n)^{\lambda}$.

Regarding the choice of the segment length $m_n$ and the number of time points $N$ used, Assumption 5(iv) and Assumption 6(iv) impose some restrictions on the rates at which these parameters have to increase to infinity. Notice that the number of time points $N$ is automatically determined by the choice of the segment length $m_n$ and of the shift $S_n$, that is, of the constant $c$. For instance, for $\lambda = 1/4$, $\delta$ is allowed to vary in the interval $(1/2, 2/3)$. Within the allowed range of possible values for $\delta$, we stress the fact that $\delta$ should be chosen on the one hand large enough in order to have within each time segment a sufficient number of observations for nonparametric estimation. On the other hand, $\delta$ should not chosen too large in order to have sufficient information on the local structure of the series, that is, a sufficient number of $N$ distinct time points. Recall that $c$ determines the degree of overlapping of the segments used and that this parameter should be kept rather small, that is, $1 \leq c \leq 3$.

## 5.2. A real-life data example

We demonstrate the capability of our method to detect changes of the spectral structure of a time series by means of a real-life data example. The data analyzed consists of the first 3072 observations of a set of tremor data recorded in the Cognitive Neuroscience Laboratory of the University of Quebec at Montreal, the purpose of the study being to compare different regions of tremor activity coming from a subject with Parkinson's disease. It has been previously analyzed by von Sachs and Neumann (2000). A plot of the first differences of the series is shown in Figure 1(a) after a Gaussian noise of standard deviation 0.01 has been added in order to break the discrete nature of the data.

We apply the test proposed using a local window of length $m_n = 256$ observations and a shift $S_n = m_n$ ($c = 1$), that is, we consider $N = 12$ time points. To smooth the global and local periodogram, the Bartlett–Priestley kernel has been used with bandwidths $h = 0.1$ and $b = 0.18$, where $h$ has been chosen as the minimizer of $CV(h)$. For this choice of parameters, we get $(m_n \sqrt{Nb} T_n - \widehat{\mu}_n)/\tau = 21.16$, which, compared with the 5% critical value of the standard Gaussian distribution, clearly leads to a rejection of the null hypothesis.

To better understand the reasons leading to the above rejection of stationarity, it is useful to consider the statistic $Q_n(u)$ defined by

$$Q_n(u) = m_n \sqrt{b} \int_{-\pi}^{\pi} V_n^2(u, \lambda) \, d\lambda, \qquad u \in (0, 1), \tag{5.1}$$

which evaluates, at any time point $u \in (0,1)$, the difference between the locally calculated sample spectral density $I_{m_n}(u, \lambda)$ and its globally estimated counterpart $\widehat{g}_h(\lambda)$. Notice that the test statistic $T_n$ is just a time-averaged version of $Q_n(u)$ calculated for $N$ different instants of time. Validity of the null hypothesis corresponds to small values of $Q_n(u)$ while large values of the same statistic pinpoint time regions where changes in the spectral structure of the series occur. Figure 1(c) shows a plot of $Q_n(t/n)$ against the rescaled time variable $t/n$. As this figure clearly indicates, the main contributions to the test statistic $T_n$ leading to the rejection of the null hypothesis are coming essentially from two different time segments. Based on the behavior of $Q_n(t_s/n)$, these segments are easily identified as follows. The first segment starts approximately at $t = 1760$ and ends at $t = 2170$, while the second starts at approximately $t = 2350$ and ends at $t = 2840$. Figure 1(b) shows the estimated time-averaged spectral density $\widehat{g}_h(\lambda)$, as well as the estimated spectral densities obtained using observations within each one of the two identified time segments. Figure 1(d) shows the estimated spectral density obtained using observations of the first part of the series. As is clearly seen from these exhibits, an increase in the variability of the series, together with a reallocation of the total power toward low-frequency components, makes the behavior of the spectral densities in the two identified time segments different from the overall behavior and from the behavior during the first part of the series.

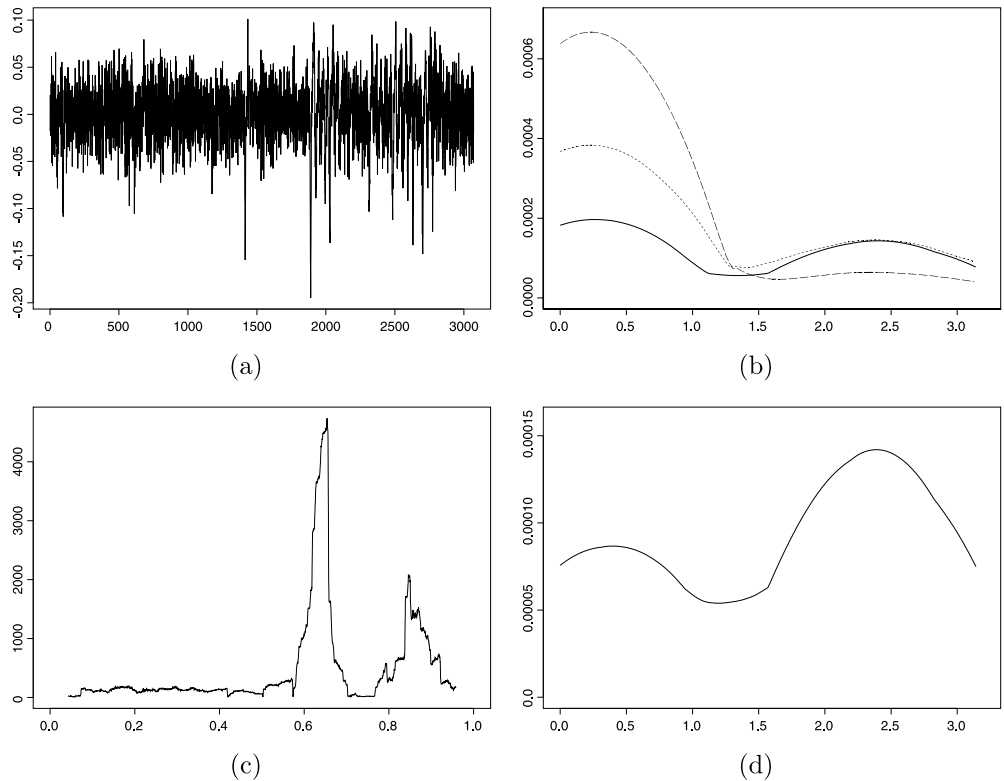

(a)   (b)

(c)   (d)

**Figure 1.** (a) Plot of the first differences of the tremor series ($n = 3071$). (b) Estimated spectral densities for different segments of the tremor series: the solid line refers to the estimated time-averaged spectral density $\widehat{g}_h(\lambda)$, the dashed line to the estimated spectral density for the time segment $t \in \{1760, 1761, \dots, 2170\}$ and the dotted line for the time segment $t \in \{2350, 2351, \dots, 2840\}$. (c) Plot of the statistic $Q_n(u_s)$ against the rescaled time parameter $u_s = t_s/3071$, $t_s = 129, 130, \dots, 2942$. (d) Estimated spectral density of the first time segment consisting of the observations $X_t$ for $t \in \{1, 2, \dots, 1759\}$.

# 6. Auxiliary results and proofs

**Lemma 6.1.** *Suppose Assumptions* 1 *and* 3 *hold. As* $n \to \infty$,

(i) $E(\widehat{g}_h(\lambda) - \int_0^1 f(u, \lambda)\,\mathrm{d}u)^2 = \mathrm{O}(n^{-1}h^{-1}) + \mathrm{O}(h^4)$, *where the* $\mathrm{O}(\cdot)$ *terms are uniformly in* $\lambda \in [-\pi, \pi]$;

(ii) $\sup_{\lambda \in [-\pi, \pi]} |\widehat{g}_h(\lambda) - \int_0^1 f(u, \lambda)\,\mathrm{d}u| \to 0$ *in probability.*

**Proof.** To prove (i), let

$$d_n(\lambda_k) = \frac{1}{\sqrt{2\pi n}} \sum_{t=1}^n X_{t,n} \exp\{-\mathrm{i}\lambda_k t\}$$

and $\Delta_m(\lambda) = \sum_{s=1}^m \exp\{-\mathrm{i}\lambda s\}$. Note that $|\Delta_m(\lambda)| \leq \pi L_m(\lambda)$, where the function $L_m(\omega)$ is the periodic extension of $L_m(\omega) = n\mathbf{1}_{\{|\omega| \leq n^{-1}\}}(\omega) + |\omega|^{-1}\mathbf{1}_{\{n^{-1} < |\omega| \leq \pi\}}(\omega)$; cf. Dahlhaus (1997). Using summation by parts, we get

$$\sum_{t=1}^n A(t/n, \lambda_k)\mathrm{e}^{-\mathrm{i}(\lambda - \lambda_k)t} = A(1, \lambda_k)\mathrm{O}(L_n(\lambda_k - \lambda)) + \mathrm{O}\Big(\sup_u |A'(u, \lambda_k)| L_n(\lambda_k - \lambda)\Big)$$

$$= \mathrm{O}(L_n(\lambda_k - \lambda)). \tag{6.1}$$

Furthermore, because

$$\left| \sum_{t=1}^n (A_{t,n}(\lambda) - A(t/n, \lambda_k))\mathrm{e}^{-\mathrm{i}(\lambda - \lambda_k)t} \right|$$

$$\leq \sum_{t=1}^n |A_{t,n}(\lambda) - A(t/n, \lambda)| + \left| \sum_{t=1}^n (A(t/n, \lambda) - A(t/n, \lambda_k))\mathrm{e}^{-\mathrm{i}(\lambda - \lambda_k)t} \right|$$

$$\leq n \sup_{t,n} |A_{t,n}(\lambda) - A(t/n, \lambda)| + \sup_u |A'(u, \lambda)||\lambda - \lambda_k||\Delta_n(\lambda_k - \lambda)|$$

$$\leq K^{-1}(1 + |\lambda_k - \lambda| L_n(\lambda_k - \lambda)),$$

it follows that

$$\left| \sum_{t=1}^{n} (A_{t,n}(\lambda) - A(t/n, \lambda_k)) \mathrm{e}^{-\mathrm{i}(\lambda - \lambda_k)t} \right| = \mathrm{O}(1). \tag{6.2}$$

By (6.1), (6.2) and the fact that $\int_{-\pi}^{\pi} L_n(\lambda - \lambda_k) \, \mathrm{d}\lambda = \mathrm{O}(\ln(n))$ (cf. Dahlhaus (1997), Lemma A.4), we get

$$\begin{aligned}
E(I_n(\lambda_k)) &= cum(d_n(\lambda_k), d_n(-\lambda_k)) \\
&= \frac{1}{2\pi n} \int_{-\pi}^{\pi} \sum_{t=1}^{n} \sum_{s=1}^{n} A_{t,n}(\lambda) \overline{A_{s,n}(\lambda)} \mathrm{e}^{-\mathrm{i}(\lambda + \lambda_k)(t-s)} \, \mathrm{d}\lambda \\
&= \frac{1}{2\pi n} \int_{-\pi}^{\pi} \sum_{t=1}^{n} \sum_{s=1}^{n} A(t/n, \lambda_k) \overline{A(s/n, \lambda_k)} \mathrm{e}^{-\mathrm{i}(\lambda + \lambda_k)(t-s)} \, \mathrm{d}\lambda + \mathrm{O}(\ln(n)n^{-1}) \\
&= \int_{0}^{1} f(u, \lambda_k) \, \mathrm{d}u + \mathrm{O}(\ln(n)n^{-1}).
\end{aligned}$$

Let $\omega = \lambda_j$ and $\lambda = \lambda_k$ and consider

$$\begin{aligned}
Cov(I_n(\lambda), I_n(\omega)) &= cum(d_n(\omega), d_n(\lambda)) cum(d_n(-\omega), d_n(-\lambda)) \\
&\quad + cum(d_n(\omega) d_n(-\lambda)) cum(d_n(-\omega), d_n(\lambda)) \\
&\quad + cum(d_n(\omega), d_n(-\omega), d_n(\lambda), d_n(-\lambda)).
\end{aligned} \tag{6.3}$$

For $\omega \neq \lambda$, we have, for the first term in (6.3),

$$\begin{aligned}
cum(d_n(\omega), d_n(\lambda)) &= \frac{1}{2\pi n} \sum_{s=1}^{n} \sum_{k=-n+s}^{n-s} Cov(X_{s,n}, X_{s+k,n}) \mathrm{e}^{-\mathrm{i}\omega k} \mathrm{e}^{-\mathrm{i}(\lambda + \omega)s} \\
&= \frac{1}{2\pi n} \sum_{s=1}^{n} \sum_{k=-n+s}^{n-s} \mathrm{c}(s/n, k) \mathrm{e}^{-\mathrm{i}\omega k} \mathrm{e}^{-\mathrm{i}(\lambda + \omega)s} + \mathrm{O}(n^{-1}) \\
&= \frac{1}{2\pi n} \left( \sum_{k=0}^{n-1} \sum_{s=1}^{n-k} c(s/n, k) \mathrm{e}^{-\mathrm{i}\lambda k} \mathrm{e}^{-\mathrm{i}(\lambda + \omega)s} \right. \\
&\quad \left. + \sum_{k=-n+1}^{-1} \sum_{s=1}^{n-|k|} c(s/n, k) \mathrm{e}^{-\mathrm{i}\lambda k} \mathrm{e}^{-\mathrm{i}(\lambda + \omega)s} \right) + \mathrm{O}(n^{-1}),
\end{aligned}$$

from which we conclude, using

$$\left| \frac{1}{n} \sum_{k=0}^{n-1} \sum_{s=1}^{n-k} c(s/n, k) \mathrm{e}^{-\mathrm{i}\lambda k} \mathrm{e}^{-\mathrm{i}(\lambda + \omega)s} \right| \leq \frac{1}{n} \sum_{k=0}^{n-1} \left| \sum_{s=1}^{n-k} \mathrm{e}^{-\mathrm{i}(\lambda + \omega)s} \right| |c(s/n, k)|$$

$$\leq \frac{C}{n} \sum_{k=0}^{n-1} |k| \sum_{j=-\infty}^{\infty} \frac{1}{l(j)l(j+k)}$$

$$\leq \frac{2C}{\sqrt{n}} \left( \sum_{k=-\infty}^{\infty} |k|^{1/2} \frac{1}{l(k)} \right) \left( \sum_{j=-\infty}^{\infty} \frac{1}{l(j)} \right) = \mathrm{O}(n^{-1/2}),$$

that $cum(d_n(\omega), d_n(\lambda)) = \mathrm{O}(n^{-1/2})$ uniformly in $\omega$ and $\lambda$. Notice that the above arguments remain valid if $\omega$ (resp. $\lambda$) is replaced by $-\omega$ (resp. $-\lambda$) or both. Furthermore,

$$|cum(d_n(\omega), d_n(-\omega), d_n(\lambda), d_n(-\lambda))|$$

$$\leq n^{-2} C \sum_{t_1=1}^{n} \sum_{t_2=1}^{n} \sum_{s_1=1}^{n} \sum_{s_2=1}^{n} \sum_{j=-\infty}^{\infty} |a_{t_1,n}(j)||a_{t_2,n}(j+t_2-t_1)|$$

$$\times |a_{s_1,n}(j+s_1-t_1)||a_{s_2,n}(j+s_2-t_1)|$$

$$\leq n^{-1} C K^4 \left( \sum_{j=-\infty}^{\infty} l^{-1}(j) \right)^4 = \mathrm{O}(n^{-1}).$$

Thus, $Cov(I_n(\lambda_k), I_n(\lambda_j)) = \mathrm{O}(n^{-1})$ for $\lambda_j \neq \lambda_k$, and for $\lambda_j = \lambda_k$, we get (along the same lines) that $Var(I_n(\lambda_j)) = \mathrm{O}(1)$ uniformly in $\lambda_j$. Using these results, the properties of the kernel $K$ (see Assumption 2) and the smoothness properties of $g$ (see Assumption 1), it follows by standard arguments that $E(\widehat{g}_h(\lambda)) = g(\lambda) + \mathrm{O}(h^2)$ and $Var(\widehat{g}_h(\lambda)) = \mathrm{O}(n^{-1}h^{-1})$ uniformly in $\lambda$.

To establish (ii), notice that

$$\left| \widehat{g}_h(\lambda) - \int_0^1 f(u, \lambda) \, \mathrm{d}u \right| \leq \left| \frac{1}{n} \sum_{j=-L_n}^{L_n} K_h(\lambda - \omega_j) \left[ I_n(\omega_j) - \frac{1}{n} \sum_{t=1}^{n} f\left( \frac{t}{n}, \omega_j \right) \right] \right|$$

$$+ \left| \frac{1}{n} \sum_{j=-L_n}^{L_n} K_h(\lambda - \omega_j) \frac{1}{n} \sum_{t=1}^{n} \left[ f\left( \frac{t}{n}, \omega_j \right) - f\left( \frac{t}{n}, \lambda \right) \right] \right|$$

$$+ \left| \frac{1}{n} \sum_{t=1}^{n} f\left( \frac{t}{n}, \lambda \right) \left( \frac{1}{n} \sum_{j=-L_n}^{L_n} K_h(\lambda - \omega_j) - 1 \right) \right| + \mathrm{O}(n^{-1}).$$

The first term on the right-hand sight goes to zero because of part (i), Markov's inequality and a standard discretization argument (cf. the proof of Theorem A1 of Franke and Härdle (1992)). Using $|f(u, \lambda_1) - f(u, \lambda_2)| \leq C|\lambda_1 - \lambda_2|$, the second term is bounded by $Cn^{-2} \sum_{j=-L_n}^{L_n} \sum_{t=1}^{n} K_h(\lambda - \omega_j) \pi h = \mathrm{O}(h)$. For the third term, we have, by Assumption 3,

$$\left| \left( \frac{1}{n} \sum_{j=-L_n}^{L_n} K_h(\lambda - \omega_j) - 1 \right) \frac{1}{n} \sum_{t=1}^{n} f\left( \frac{t}{n}, \lambda \right) \right| \leq \frac{2\pi C}{n^2 h} \sum_{t=1}^{n} f\left( \frac{t}{n}, \lambda \right) = \mathrm{O}(n^{-1}h^{-1}). \qquad \square$$

**Lemma 6.2.** *Let Assumption* 1 *be true and* $E|\varepsilon_1|^k < \infty$ *for* $k \in \mathbb{N}$. *Then*

$$cum(d_n(\lambda_1), d_n(\lambda_2), \ldots, d_n(\lambda_k)) = \frac{(2\pi)^{k-1}}{n^{k/2}} \sum_{t=1}^{n} f_k\left(\frac{t}{n}, \lambda_1, \lambda_2, \ldots, \lambda_{k-1}\right) \mathrm{e}^{-\mathrm{i}t\sum_{s=1}^{k}\lambda_s}$$

$$+ \, \mathrm{O}(\ln^{k-1}(n)n^{-k/2}),$$

*where the error term is uniform in* $\lambda_1, \lambda_2, \ldots, \lambda_k$ *and for* $u \in [0, 1]$,

$$f_k(u, \lambda_1, \lambda_2, \ldots, \lambda_{k-1}) = \frac{cum_k(\varepsilon_t)}{(2\pi)^{k-1}} A(u, \lambda_1) A(u, \lambda_2) \cdots A(u, \lambda_{k-1}) A\left(u, \sum_{s=1}^{k-1} \lambda_s\right)$$

*is the kth order local cumulant spectrum of the process* $\{\widetilde{X}_t(u), t \in \mathbb{Z}\}$, *where* $\widetilde{X}_t(u) = \sum_{j=-\infty}^{\infty} a(u, j)\varepsilon_{t-j}$ *and* $A(u, \lambda) = \sum_{j=-\infty}^{\infty} a(u, j)\exp\{-\mathrm{i}j\lambda\}$.

**Proof.** We use

$$cum(d_n(\lambda_1), d_n(\lambda_2), \ldots, d_n(\lambda_k))$$

$$= \frac{cum_k(\varepsilon_t)}{2\pi^{k-1}n^{k/2}} \sum_{t_1,\ldots,t_k} \int \prod_{s=1}^{k-1} A_{t_s,n}(\omega_s) \mathrm{e}^{-\mathrm{i}(\lambda_s - \omega_s)t_s}$$

$$\times A_{t_k,n}\left(-\sum_{s=1}^{k-1} \omega_s\right) \mathrm{e}^{-\mathrm{i}(\lambda_k + \sum_{s=1}^{k-1} \omega_s)t_k} \, \mathrm{d}\omega_1 \cdots \mathrm{d}\omega_{k-1}$$

and successively replace $A_{t_s,n}(\omega_s)$ by $A(t_s/n, \lambda_s)$ for $s = 1, 2, \ldots, k-1$ and $A_{t_k,n}(-\sum_{s=1}^{k-1}\omega_s)$ by $A(t_s/n, -\sum_{s=1}^{k-1}\lambda_s)$. We then get that the expression on the right-hand side above equals

$$\frac{cum_k(\varepsilon_t)}{2\pi^{k-1}n^{k/2}} \sum_{t_1,\ldots,t_k} \int \prod_{s=1}^{k-1} A(t_s/n, \lambda_s) \mathrm{e}^{-\mathrm{i}(\lambda_s - \omega_s)t_s}$$

$$\times A\left(t_k/n, -\sum_{s=1}^{k-1} \lambda_s\right) \mathrm{e}^{-\mathrm{i}(\lambda_k + \sum_{s=1}^{k-1} \omega_s)t_k} \, \mathrm{d}\omega_1 \cdots \mathrm{d}\omega_{k-1}$$

$$= \frac{cum_k(\varepsilon_t)}{n^{k/2}} \sum_{t=1}^{n} \prod_{s=1}^{k-1} A(t/n, \lambda_s) A\left(t/n, -\sum_{s=1}^{k-1} \lambda_s\right) \mathrm{e}^{-\mathrm{i}t\sum_{s=1}^{k}\lambda_s},$$

plus a remainder $R_n$ which, using bounds similar to (6.1) and (6.2) and Lemma A.4 of Dahlhaus (1997), is shown to be, after tedious but straightforward calculations, bounded by

$$Const. \frac{cum_k(\varepsilon_t)}{2\pi^{k-1}n^{k/2}} \int \prod_{s=1}^{k-1} L_n(\lambda_s - \omega_s) \, \mathrm{d}\omega_1 \, \mathrm{d}\omega_2 \cdots \mathrm{d}\omega_{k-1} = \mathrm{O}(\ln^{k-1}(n)/n^{k/2}).$$

□

We next show that under the assumptions of Theorem 3.1, replacing $\widehat{g}_h(\lambda)$ by $g(\lambda) = f(\lambda)$ does not affect the asymptotic distribution of $T_n$. For this, write under the null hypothesis $\widehat{f}_h(\lambda)$ for $\widehat{g}_h(\lambda)$ and let $L_n(u,\lambda) = (I_{m_n}(u,\lambda)/f(\lambda) - 1)$ and $G_n(\lambda) = (\widehat{f}_h(\lambda)/f(\lambda) - 1)$. Further, let

$$\widetilde{T}_n = \frac{1}{N} \sum_{s=1}^{N} \int_{-\pi}^{\pi} \widetilde{V}_n^2\left(\frac{t_s}{n}, \lambda\right) d\lambda, \qquad (6.4)$$

where $\widetilde{V}_n(u,\lambda)$ is defined in (2.11). We then have the following lemma.

**Lemma 6.3.** *Suppose the assumptions of Theorem 3.1 hold. As $n \to \infty$,*

$$T_n - \widetilde{T}_n = o_P(m_n^{-1} N^{-1/2} b^{-1/2}).$$

**Proof.** Using a Taylor series argument and Lemma 6.1, we get

$$m_n \sqrt{bN} T_n$$

$$= m_n \sqrt{bN} \widetilde{T}_n - 2 \frac{m_n \sqrt{b}}{\sqrt{N}} \sum_{s=1}^{N} \int_{-\pi}^{\pi} \left\{ \frac{1}{m_n} \sum_j K_b(\lambda - \lambda_j) L_n(u_s, \lambda_j) \right\}$$

$$\times \left\{ \frac{1}{m_n} \sum_l K_b(\lambda - \lambda_l)(L_n(u_s, \lambda_l) + 1) G_n(\lambda_l) \right\} d\lambda$$

$$+ O_P\left(\frac{m_n \sqrt{bN}}{nh}\right) + O_P(m_n \sqrt{bN} h^4),$$

where the first $O_P(\cdot)$ term goes to zero for $\delta < (1-2\lambda)/(1-\lambda)$ and the second for $\delta < (8\lambda - 1)/(1-\lambda)$. Denoting by $-2R_n$ the second term on the right-hand side of the above equation, we verify that

$$R_n = \frac{\sqrt{b}}{m_n \sqrt{N}} \sum_{s=1}^{N} \int_{-\pi}^{\pi} \sum_j \sum_l K_b(\lambda - \lambda_j) K_b(\lambda - \lambda_l) L_n(u_s, \lambda_j) G_n(\lambda_l) d\lambda + O_P\left(\sqrt{\frac{N}{nbh}}\right)$$

$$= \widetilde{R}_n + O_P\left(\sqrt{\frac{N}{nbh}}\right)$$

with an obvious definition for $\widetilde{R}_n$ and where $O_P(\sqrt{N/nbh}) \to 0$ for $\delta > \lambda/(1-\lambda)$. Evaluating the expectation $E(L_n(u_1, \lambda_1) G_n(\lambda_2) L_n(u_2, \lambda_3) G_n(\lambda_4))$, we can decompose $E(\widetilde{R}_n)^2$ into a sum of several terms, two typical of which are

$$E_{1,n} = \frac{b}{m_n^2 N} \sum_{s_1} \sum_{s_2} \int \int \sum_{j_1} \sum_{l_1} \sum_{j_2} \sum_{l_2} K_b(\lambda_1 - \lambda_{j_1}) K_b(\lambda_1 - \lambda_{l_1}) K_b(\lambda_2 - \lambda_{j_2})$$

$$\times K_b(\lambda_2 - \lambda_{l_2})\, d\lambda_1\, d\lambda_2$$

$$\times E[L_n(u_{s_1}, \lambda_{j_1}) L_n(u_{s_2}, \lambda_{j_2})] E[G_n(\lambda_{l_1}) G_n(\lambda_{l_2})]$$

and

$$E_{2,n} = \frac{b}{m_n^2 N} \sum_{s_1} \sum_{s_2} \int \int \sum_{j_1} \sum_{l_1} \sum_{j_2} \sum_{l_2} K_b(\lambda_1 - \lambda_{j_1}) K_b(\lambda_1 - \lambda_{l_1}) K_b(\lambda_2 - \lambda_{j_2})$$

$$\times K_b(\lambda_2 - \lambda_{l_2})\, d\lambda_1\, d\lambda_2$$

$$\times E[L_n(u_{s_1}, \lambda_{j_1}) G_n(\lambda_{l_1})] E[L_n(u_{s_2}, \lambda_{j_2}) G_n(\lambda_{l_2})].$$

Using $E(I_{m_n}(u, \lambda_j)) = f(\lambda_j) + O(m_n^{-2})$, we get

$$E_{1,n} = \frac{b}{m_n^2 N} \sum_{s_1} \sum_{s_2} \int \int \sum_{j_1} \sum_{l_1} \sum_{j_2} \sum_{l_2} K_b(\lambda_1 - \lambda_{j_1}) K_b(\lambda_1 - \lambda_{l_1}) K_b(\lambda_2 - \lambda_{j_2})$$

$$\times K_b(\lambda_2 - \lambda_{l_2})\, d\lambda_1\, d\lambda_2$$

$$\times Cov(I_{m_n}(u_{s_1}, \lambda_{j_1})/f(\lambda_{j_1}), I_{m_n}(u_{s_2}, \lambda_{j_2})/f(\lambda_{j_2}))$$

$$\times E[G_n(\lambda_{l_1}) G_n(\lambda_{l_2})] + O(Nb(nh)^{-1}).$$

Now, standard calculations (cf. Theorem 4.3.2 of Brillinger (1981)) yield

$$Cov(I_{m_n}(u_{s_1}, \lambda_j), I_{m_n}(u_{s_2}, \lambda_l))$$

$$= \frac{1}{4\pi^2 H_{2,m_n}^2(0)} f(\lambda_j) f(\lambda_l) \tag{6.5}$$

$$\times \Bigg[ \int \int H_{m_n}(\omega_1 - \lambda_j) H_{m_n}(-\omega_1 + \lambda_l) H_{m_n}(\omega_2 + \lambda_j) H_{m_n}(-\omega_2 - \lambda_l)$$

$$\times e^{iS(s_1 - s_2)(\omega_1 + \omega_2)}\, d\omega_1\, d\omega_2$$

$$+ \int \int H_{m_n}(\omega_1 - \lambda_j) H_{m_n}(-\omega_1 - \lambda_l) H_{m_n}(\omega_2 + \lambda_j) H_{m_n}(-\omega_2 + \lambda_l)$$

$$\times e^{iS(s_1 - s_2)(\omega_1 + \omega_2)}\, d\omega_1\, d\omega_2$$

$$+ \int \int \int f_4(\omega_1, \omega_2, \omega_3) e^{iS(s_1 - s_2)(\omega_2 + \omega_3)} H_{m_n}(\omega_1 + \omega_2 + \omega_3 + \lambda_j)$$

$$\times H_{m_n}(-\omega_1 - \lambda_j) H_{m_n}(-\omega_2 + \lambda_l) H_{m_n}(-\omega_3 - \lambda_l)\, d\omega_1\, d\omega_2\, d\omega_3 \Bigg]$$

$$+ O(\log^2(m_n) m_n^{-2})$$

$$= \sum_{r=1}^{3} C_{r,n}(u_{s_1}, u_{s_2}, \lambda_j, \lambda_l) + O(\log^2(m_n) m_n^{-2}),$$

with an obvious definition for $C_{r,n}(u_{s_1}, u_{s_2}, \lambda_j, \lambda_l)$, $r = 1, 2, 3$ and where $f_4(\lambda_1, \lambda_2, \lambda_3) = (2\pi)^{-3}\kappa_4 A(e^{-i\lambda_1})A(e^{-i\lambda_2})A(e^{-i\lambda_3})A(e^{i(\lambda_1+\lambda_2+\lambda_3)})$, with $A(z) = \sum_j a(j)z^j$, denotes the fourth order cumulant spectrum, and $H_{m_n}(\lambda) = H_{1,m_n}(\lambda)$. Substituting (6.5) into $E_{1,n}$, we get by straightforward calculation that $E_{1,n} = O(m_n b(nh)^{-1} + m_n h^4 b)$. To demonstrate this, consider, for instance, the term of $E_{1,n}$ corresponding to $C_{1,n}(u_{s_1}, u_{s_2}, \lambda_{j_1}, \lambda_{j_2})$. Using

$$
\begin{aligned}
&\frac{1}{2\pi}\int_{-\pi}^{\pi} H_{m_n}(x-y_1)H_{m_n}(-x+y_2)e^{iS(s_1-s_2)x}\,\mathrm{d}x \\
&= \sum_{\substack{t_1=1 \\ t_1-t_2=S(s_1-s_2)}}^{m_n}\sum_{t_2=1}^{m_n} h\left(\frac{t_1}{m_n}\right)h\left(\frac{t_2}{m_n}\right)e^{i(y_1t_1-y_2t_2)}
\end{aligned}
\tag{6.6}
$$

and $E[G_n(\lambda_{l_1})G_n(\lambda_{l_2})] = O(n^{-1}h^{-1} + h^4)$, by Lemma 6.1, we get that this term can be bounded by

$$
\mathrm{O}\!\left(m_n^2\left(\frac{1}{nh}+h^4\right)\right)\frac{b}{m_n^2 N}
$$

$$
\times \sum_{s_1=1}^{N}\sum_{s_2=1}^{N}\int\int\sum_{j_1}\sum_{j_2}K_b(\lambda_1-\lambda_{j_1})K_b(\lambda_2-\lambda_{j_2})
$$

$$
\times \frac{1}{H_{2,m_n}^2(0)}\left|\sum_{\substack{t_1=1 \\ t_1-t_2=S(s_1-s_2)}}^{m_n}\sum_{t_2=1}^{m_n} h\left(\frac{t_1}{m_n}\right)h\left(\frac{t_2}{m_n}\right)e^{i(\lambda_{j_1}t_1-\lambda_{j_2}t_2)}\right|^2\mathrm{d}\lambda_1\,\mathrm{d}\lambda_2
$$

$$
= \mathrm{O}\!\left(m_n^2\left(\frac{1}{nh}+h^4\right)\right)\frac{b}{H_{2,m_n}^2(0)m_n^2}
$$

$$
\times \sum_{|r|<c}\left(1-\frac{|r|}{N}\right)\int\sum_{j_1,j_2}K_b(\lambda_1-\lambda_{j_1})K_b(\lambda_2-\lambda_{j_2})
$$

$$
\times \left|\sum_{t=1+rS\mathbf{1}(0\le r<c)}^{m_n+rS\mathbf{1}(-c<r<0)} h\left(\frac{t}{m_n}\right)h\left(\frac{t+rS}{m_n}\right)e^{i(\lambda_{j_1}-\lambda_{j_2})t}\right|^2\mathrm{d}\lambda_1\,\mathrm{d}\lambda_2.
$$

Approximating the sum over the Fourier frequencies by an integral and using the Lipschitz continuity of the kernel $K$, we get that the above term is bounded by

$$
\mathrm{O}\!\left(m_n^2\left(\frac{1}{nh}+h^4\right)\right)\frac{b}{H_{2,m_n}^2(0)}
$$

$$
\times \sum_{|r|<c}\int K_b(\lambda_1-x)K_b(\lambda_2-x)\,\mathrm{d}\lambda_1\,\mathrm{d}\lambda_2
$$

$$\times \left| \sum_{t=1+rS\mathbf{1}(0\le r<c)}^{m_n+rS\mathbf{1}(-c<r<0)} h\left(\frac{t}{m_n}\right) h\left(\frac{t+rS}{m_n}\right) \mathrm{e}^{\mathrm{i}(x-y)t} \right|^2 \mathrm{d}x\,\mathrm{d}y + RM_n$$

$$= \mathrm{O}\left(m_n^2\left(\frac{1}{nh}+h^4\right)\right)\mathrm{O}(bH_{2,m_n}^{-1}(0)) \int K_b(\lambda_1-x)K_b(\lambda_2-x)\,\mathrm{d}x\,\mathrm{d}\lambda_1\,\mathrm{d}\lambda_2 + RM_n$$

$$= \mathrm{O}\left(\frac{m_nb}{nh}+m_nh^4b\right) + RM_n,$$

where the remainder $RM_n$ satisfies

$$|RM_n| \le \mathrm{O}\left(m_n^2\left(\frac{1}{nh}+h^4\right)\right)\frac{1}{H_{2,m_n}^2(0)}$$

$$\times \int K_b(\lambda_1-x)|y-x|\,\mathrm{d}\lambda_1$$

$$\times \left| \sum_{t=1+rS\mathbf{1}(0\le r<c)}^{m_n+rS\mathbf{1}(-c<r<0)} h\left(\frac{t}{m_n}\right) h\left(\frac{t+rS}{m_n}\right) \mathrm{e}^{\mathrm{i}(x-y)t} \right|^2 \mathrm{d}x\,\mathrm{d}y$$

$$\le \mathrm{O}\left(m_n^2\left(\frac{1}{nh}+h^4\right)\right)\frac{1}{H_{2,m_n}^2(0)} \int\int |y-x|L_{m_n}^2(x-y)\,\mathrm{d}x\,\mathrm{d}y$$

$$= \mathrm{O}(n^{-1}h^{-1}+h^4).$$

The function $L_m(\omega)$ appearing above is the periodic extension of $L_m(\omega) = n\mathbf{1}_{\{|\omega|\le n^{-1}\}}(\omega) + |\omega|^{-1}\mathbf{1}_{\{n^{-1}<|\omega|\le\pi\}}(\omega)$ and the last equality follows using Lemma A.4 of Dahlhaus ([1997]). To bound $E_{2,n}$, we use

$$E[L_n(u_s,\lambda_j)G_n(\lambda_l)] = f^{-1}(\lambda_j)f^{-1}(\lambda_l)\,cum(I_{m_n}(u_s,\lambda_j),\widehat{f}(\lambda_l)) + \mathrm{O}(h^2m_n^{-2})$$

and $cum(I_{m_n}(u_s,\lambda_j),I_n(\lambda_r)) = \mathrm{O}(m_n^{-1/2}n^{-1/2})$ for $\lambda_j \ne \lambda_r$, where the $\mathrm{O}(\cdot)$ term is uniformly in $u_s$, $\lambda_j$ and $\lambda_r$. We then have

$$\frac{\sqrt{b}}{m_n n\sqrt{N}} \sum_s \int \sum_j \sum_l \sum_r K_b(\lambda-\lambda_j)K_b(\lambda-\lambda_l)K_h(\lambda_l-\lambda_r)\,\mathrm{d}\lambda\,\frac{1}{f(\lambda_j)f(\lambda_l)}$$

$$\times cum(I_{m_n}(u_s,\lambda_j),I_n(\lambda_r))$$

$$= \mathrm{O}(m_n^{1/2}b^{1/2}N^{1/2}n^{-1/2}) = \mathrm{O}(b^{1/2})$$

and, therefore, $E_{2,n} = \mathrm{O}(b)$. $\qquad\square$

**Proof of Theorem 3.1.** It suffices to show that $m_n\sqrt{Nb}\widetilde{T}_n - \mu_n$ converges to the desired Gaussian distribution. Using (6.5) and $E(I_{m_n}(u_s,\lambda_j)) = f(\lambda_j) + \mathrm{O}(m_n^{-2})$, we get

$$m_n \sqrt{N} b E(\widetilde{T}_n)$$

$$= \frac{\sqrt{b}}{m_n \sqrt{N}} \sum_{s=1}^{N} \sum_{j_1} \sum_{j_2} \int_{-\pi}^{\pi} K_b(\lambda - \lambda_{j_1}) K_b(\lambda - \lambda_{j_2}) \, d\lambda \frac{1}{f(\lambda_{j_1}) f(\lambda_{j_2})}$$

$$\times Cov(I_{m_n}(u_s, \lambda_{j_1}), I_{m_n}(u_s, \lambda_{j_2})) + O(N^{1/2} m_n^{-3/2})$$

$$= \frac{\sqrt{bN}}{m_n H_{2,m_n}^2(0)} \sum_{j_1} \sum_{j_2} \int_{-\pi}^{\pi} K_b(\lambda - \lambda_{j_1}) K_b(\lambda - \lambda_{j_2}) \, d\lambda$$

$$\times \{|H_{2,m_n}(\lambda_{j_1} + \lambda_{j_2})|^2 + |H_{2,m_n}(\lambda_{j_1} - \lambda_{j_2})|^2 + H_{4,m_n}(0)(\eta - 3)\}$$

$$+ O(N^{1/2} m_n^{-3/2})$$

$$= E_{1,n} + E_{2,n} + E_{3,n} + O(N^{1/2} m_n^{-3/2})$$

with an obvious definition for the $E_{i,n}$'s. For $E_{1,n}$, we have

$$E_{1,n} = \frac{m_n \sqrt{bN}}{4\pi^2 H_{2,m_n}^2(0)} \int \int \int K_b(v + 2x) K_b(v) |H_{2,m_n}(\omega)|^2 \, d\omega \, dv \, dx + O(N^{1/2} b^{-1/2} m_n^{-1}),$$

where $O(N^{1/2} b^{-1/2} m_n^{-1}) \to 0$ for $\delta > 1/(3 - \lambda)$ and

$$\left| \frac{m_n \sqrt{bN}}{4\pi^2 H_{2,m_n}^2(0)} \int K_b(v + 2x) K_b(v) |H_{2,m_n}(\omega)|^2 \, d\omega \, dv \, dx - \frac{H_4 \sqrt{bN}}{4\pi H_2^2} \int (K * K)(y) \, dy \right| \to 0.$$

Notice that the $O(N^{1/2} b^{-1/2} m_n^{-1})$ above appears because by the Lipschitz properties of $K$, the fact that $|H_{k,m_n}(\lambda)| \le C L_{m_n}(\lambda)$, where the function $L_{m_n}(\lambda)$ has been introduced in the proof of Lemma 6.1, and by Lemma A.4 of Dahlhaus (1997), we get

$$\frac{m_n \sqrt{bN}}{4\pi^2 H_{2,m_n}^2(0)} \int K_b(v)(K_b(v + 2x + \omega) - K_b(v + 2x)) |H_{2,m_n}(\omega)|^2 \, d\omega \, dv \, dx$$

$$= O(N^{1/2} b^{-1/2} m_n^{-1}).$$

To evaluate $E_{2,n}$ we first write

$$E_{2,n} = \frac{m_n \sqrt{bN}}{4\pi^2 H_{2,m_n}^2(0)} \left\{ \int \int \int K_b(v)(K_b(v + \omega) - K_b(v)) |H_{2,m_n}(\omega)|^2 \, dv \, dx \, d\omega \right.$$

$$\left. + 2\pi \int \int K_b^2(v) |H_{2,m_n}(\omega)|^2 \, dv \, d\omega \right\} + O(N^{1/2} b^{-1/2} m_n^{-1}).$$

As for $E_{1,n}$, for the first term on the right-hand side, we get the bound

$$C \frac{m_n \sqrt{N}}{2\pi H_{2,m_n}^2(0) \sqrt{b}} \int K_b(v) \, dv \int |\omega| |H_{2,m_n}(\omega)|^2 \, d\omega = O(N^{1/2} b^{-1/2} m_n^{-1}),$$

while for the second term, we obtain

$$\left| \frac{m_n\sqrt{N}}{2\pi H_{2,m_n}^2(0)b^{3/2}} \int K^2(v/b)\,\mathrm{d}v \int |H_{2,m_n}(\omega)|^2\,\mathrm{d}\omega - \frac{H_4\sqrt{N}}{H_2^2\sqrt{b}} \int K^2(x)\,\mathrm{d}x \right| \to 0.$$

The proof that $|m_n\sqrt{bN}E(\widetilde{T}_n) - \mu_n| \to 0$ concludes at this point because

$$\left| \frac{2\pi H_{4,m_n}(0)\sqrt{bN}}{m_n H_{2,m_n}^2(0)}(\eta - 3) - \frac{2\pi H_4\sqrt{bN}}{H_2^2}(\eta - 3) \right| \to 0.$$

To obtain the variance of $m_n\sqrt{bN}\widetilde{T}_n$, notice that

$$Var(\widetilde{T}_n) = \frac{1}{m_n^4 N^2}$$

$$\times \sum_{s_1,s_2}\sum_{j_1,j_2}\sum_{l_1,l_2} \int\int K_b(\lambda_1 - \lambda_{j_1})K_b(\lambda_1 - \lambda_{j_2})K_b(\lambda_2 - \lambda_{l_1})K_b(\lambda_2 - \lambda_{l_2})\,\mathrm{d}\lambda_d\lambda_2$$

$$\times \{E[L_n(u_{s_1}, \lambda_{j_1})L_n(u_{s_2}, \lambda_{l_1})]E[L_n(u_{s_1}, \lambda_{j_2})L_n(u_{s_2}, \lambda_{l_2})]$$

$$+ E[L_n(u_{s_1}, \lambda_{j_1})L_n(u_{s_2}, \lambda_{l_2})]E[L_n(u_{s_1}, \lambda_{j_2})L_n(u_{s_2}, \lambda_{l_1})] \quad (6.7)$$

$$+ cum[L_n(u_{s_1}, \lambda_{j_1}), L_n(u_{s_1}, \lambda_{j_2}), L_n(u_{s_2}, \lambda_{l_1}), L_n(u_{s_2}, \lambda_{l_2})]\}$$

$$+ \mathrm{o}(m_n^{-2}N^{-1}b^{-1})$$

$$= V_{1,n} + V_{2,n} + V_{3,n} + \mathrm{o}(m_n^{-2}N^{-1}b^{-1})$$

with an obvious definition for $V_{i,n}$, $i = 1, 2, 3$. Using $L_n(u, \lambda_j) = f^{-1}(\lambda_j)(I_{m_n}(u, \lambda_j) - E(I_{m_n}(u, \lambda_j))) + \mathrm{O}(m_n^{-2})$, we get

$$V_{1,n} = \frac{1}{N^2}\sum_{s_1=1}^{N}\sum_{s_2=1}^{N}\int\int\left[\frac{1}{m_n^2}\sum_j\sum_l K_b(\lambda_1 - \lambda_j)K_b(\lambda_2 - \lambda_l)f^{-1}(\lambda_j)f^{-1}(\lambda_l)\right.$$

$$\left. \times Cov(I_{m_n}(u_{s_1}, \lambda_j), I_{m_n}(u_{s_2}, \lambda_l))\right]^2\mathrm{d}\lambda_1\,\mathrm{d}\lambda_2 \quad (6.8)$$

$$+ \mathrm{O}(m_n^{-5}b^{-1}).$$

Substituting (6.5) into (6.8) and ignoring the $\mathrm{O}(\log^2(m_n)m_n^{-2})$ term yields

$$V_{1,n} = \frac{1}{N^2 m_n^4}\sum_{q=1}^{3}\sum_{s_1,s_2}\int\int\left[\sum_j\sum_l K_b(\lambda_1 - \lambda_j)K_b(\lambda_2 - \lambda_l)C_{q,n}(u_{s_1}, u_{s_2}, \lambda_j, \lambda_l)\right]^2\mathrm{d}\lambda_1\,\mathrm{d}\lambda_2$$

$$+ \frac{2}{N^2 m_n^4}\sum_{q_1=1}^{3}\sum_{q_2=q_1+1}^{3}\sum_{s_1,s_2}\int\int\left[\sum_j\sum_l K_b(\lambda_1 - \lambda_j)K_b(\lambda_2 - \lambda_l)\right.$$

$$\times C_{q_1,n}(u_{s_1}, u_{s_2}, \lambda_j, \lambda_l) \Big]$$

$$\times \Big[ \sum_j \sum_l K_b(\lambda_1 - \lambda_j) K_b(\lambda_2 - \lambda_l)$$

$$\times C_{q_2,n}(u_{s_1}, u_{s_2}, \lambda_j, \lambda_l) \Big] d\lambda_1 d\lambda_2$$

$$+ o(m_n^{-2} N^{-1} b^{-1})$$

$$= \sum_{l=1}^{6} V_{1,n}^{(l)} + o(m_n^{-2} N^{-1} b^{-1})$$

with an obvious definition for $V_{1,n}^{(l)}$, $l = 1, 2, \ldots, 6$. Using (6.6), we obtain

$$
V_{1,n}^{(1)} = \frac{1}{N H_{2,m_n}^4(0)}
$$

$$
\times \sum_{|r| < c} (1 - |r|/N)
$$

$$
\times \int\int \Big[ \frac{1}{m_n^2} \sum_j \sum_l K_b(\lambda_1 - \lambda_j) K_b(\lambda_2 - \lambda_l) \tag{6.9}
$$

$$
\times \sum_{t_1, t_2 = 1 + rS\mathbf{1}(0 \le r < c)}^{m_n + rS\mathbf{1}(-c < r < 0)} h\left(\frac{t_1}{m_n}\right) h\left(\frac{t_1 - rS}{m_n}\right)
$$

$$
\times h\left(\frac{t_2}{m_n}\right) h\left(\frac{t_2 - rS}{m_n}\right)
$$

$$
\times e^{-i(\lambda_j - \lambda_l)(t_1 - t_2)} \Big]^2 d\lambda_1 d\lambda_2
$$

which, multiplied with $m_n^2 N b$, converges, as $n \to \infty$, to

$$
m_n^2 N b V_{1,n}^{(1)} \to \frac{1}{4\pi^2 H_2^4} \sum_{|s| < c} \left( \int_{[0, 1 - |s|/c]} h^2(x) h^2(x + |s|/c) \, dx \right)^2 \int_{[-2\pi, 2\pi]} (K * K)^2(y) \, dy.
$$

A repetition of the above arguments yields that $m_n^2 N b V_{1,n}^{(2)}$ converge to the same limit. For the term $V_{1,n}^{(3)}$, notice first that using $|H_{m_n}(\omega_1 + \omega_2 + \omega_3 + \lambda_j) - H_{m_n}(\omega_1 + \lambda_j)| \le$

$C|\omega_2 + \omega_3|$ and Lemma A.4 of Dahlhaus ([1997](#)), we get that

$$\frac{1}{H_{2,m_n}^2(0)m_n^2} \int \int \int \sum_j \sum_l K_b(\lambda_1 - \lambda_j)K_b(\lambda_1 - \lambda_l)|H_{m_n}(\omega_1 + \lambda_j)|^2$$

$$\times |H_{m_n}(\omega_2 + \lambda_l)||H_{m_n}(-\omega_2 - \lambda_l)| \, \mathrm{d}\omega_1 \, \mathrm{d}\omega_2 \, \mathrm{d}\omega_3 = \mathrm{O}(H_{2,m_n}^{-1}(0)).$$

This, together with straightforward calculations, yields $m_n^2 NbV_{1,n}^{(3)} = \mathrm{O}(b)$. For the cross-product terms, we have

$$m_n^2 NbV_{1,n}^{(4)}$$

$$= \frac{2m_n^2 b}{H_{2,m_n}^4(0)}$$

$$\times \sum_{|r| < c} (1 - |r|/N) \int \int \left[ \frac{1}{m_n^2} \sum_{j_1} \sum_{l_1} K_b(\lambda_1 - \lambda_{j_1}) K_b(\lambda_2 - \lambda_{l_1}) \right.$$

$$\times \left| \sum_{t=1+rS\mathbf{1}(0 \le r < c)}^{m_n+rS\mathbf{1}(-c < r < 0)} h\left(\frac{t}{m_n}\right) h\left(\frac{t - rS}{m_n}\right) \right.$$

$$\left. \left. \times \mathrm{e}^{-\mathrm{i}(\lambda_{j_1} - \lambda_{l_1})t_1} \right|^2 \right] \mathrm{d}\lambda_1 \, \mathrm{d}\lambda_2$$

$$\to \frac{1}{2\pi^2 H_2^4} \sum_{|s| < c} \left( \int_{[0,1-|s|/c]} h^2(x)h^2(x + |s|/c) \, \mathrm{d}x \right)^2 \int_{[-2\pi, 2\pi]} (K * K)^2(y) \, \mathrm{d}y.$$

Similar arguments show that $m_n^2 NbV_{1,n}^{(5)} = \mathrm{O}(b)$ and $m_n^2 NbV_{1,n}^{(6)} = \mathrm{O}(b)$.

Combining the previous results for $V_{1,n}^{(l)}$, $l = 1, 2, \ldots, 6$, we get that

$$m_n^2 NbV_{1,n} \to \frac{1}{\pi^2 H_2^4} \sum_{|s| < c} \left( \int_{[0,1-|s|/c]} h^2(x)h^2(x + |s|/c) \, \mathrm{d}x \right)^2$$

$$\times \int_{[-2\pi, 2\pi]} (K * K)^2(y) \, \mathrm{d}y \tag{6.10}$$

as $n \to \infty$. Since the structure of the term $V_{2,n}$ given in decomposition (6.7) is identical to that of $V_{1,n}$, $m_n^2 NbV_{2,n}$ converges to the same limit as (6.10). The proof that $Var(m_n\sqrt{Nb}\widetilde{T}_n)) \to \tau^2$ is then completed because $m_n^2 NbV_{3,n} \to 0$ as $n \to \infty$; see Paparoditis ([2006](#)) for details.

To show convergence of $m_n\sqrt{Nb}(T_n - ET_n)$ to the desired Gaussian distribution, notice first that as in Theorem 10.3.1 of Brockwell and Davis ([1991](#)), we have $I_{m_n}(u, \lambda_j)/f(\lambda_j) = I_\varepsilon(u, \lambda_j)/\sigma^2 + R_{m_n}(u, \lambda_j)$, where $\sup_j E(R_{m_n}(u, \lambda_j))^2 = \mathrm{O}(m_n^{-1})$, $R_{m_n}(u, \lambda_j) = a(\mathrm{e}^{-\mathrm{i}\lambda}) \times$ $J_\varepsilon(u,\lambda)Y_{m_n}(u,\lambda) + a(\mathrm{e}^{\mathrm{i}\lambda})J_\varepsilon(u,-\lambda)Y_{m_n}(u,-\lambda) + |Y_{m_n}(u,\lambda)|^2$, $J_\varepsilon(u,\lambda) = (2\pi H_{2,m_n}(0))^{-1/2} \times$
$\sum_{t=1}^{m_n} h_{t,m_n}\varepsilon_{t+[un]-M_n-1}\mathrm{e}^{-\mathrm{i}\lambda t}$,

$$Y_{m_n}(u,\lambda) = \frac{1}{\sqrt{2\pi H_{2,m_n}(0)}} \sum_{j=-\infty}^{\infty} a(j)\mathrm{e}^{-\mathrm{i}\lambda j}U_{m_n,j}(u,\lambda)$$

and

$$U_{m_n,j}(u,\lambda) = \sum_{t=1}^{m_n} h_{t,m_n}(\varepsilon_{t+[un]-M_n-1-j}\mathrm{e}^{\mathrm{i}\lambda(t-j)} - \varepsilon_{t+[un]-M_n-1}\mathrm{e}^{-\mathrm{i}\lambda t})$$

$$=: U_{m_n,j}^{(1)}(u,\lambda) + U_{m_n,j}^{(2)}(u,\lambda). \qquad (6.11)$$

Let

$$L_{n,\varepsilon}(u,\lambda) = m_n^{-1}\sum_j K_b(\lambda - \lambda_j)(I_\varepsilon(u,\lambda_j) - 1)$$

and

$$W_n(u,\lambda) = m_n^{-1}\sum_j K_b(\lambda - \lambda_j)R_{m_n}(u,\lambda_j).$$

We then have $m_n\sqrt{Nb}(T_n - ET_n) = N^{-1/2}\sum_{s=1}^N Z_{s,n} + R_{1,n} + 2R_{2,n}$, where

$$Z_{s,n} = m_n\sqrt{b}\int (L_{n,\varepsilon}^2(u_s,\lambda) - EL_{n,\varepsilon}^2(u_s,\lambda))\,\mathrm{d}\lambda,$$

$$R_{1,n} = \frac{m_n\sqrt{b}}{\sqrt{N}}\sum_{s=1}^N \int (W_n^2(u_s,\lambda) - EW_n^2(u_s,\lambda))\,\mathrm{d}\lambda$$

and

$$R_{2,n} = \frac{m_n\sqrt{b}}{\sqrt{N}}\sum_{s=1}^N \int (W_n(u_s,\lambda)L_{n,\varepsilon}(u_s,\lambda) - E(W_n(u_s,\lambda)L_{n,\varepsilon}(u_s,\lambda)))\,\mathrm{d}\lambda.$$

Observe that $\{Z_{s,n}, s = 1, 2, \ldots\}$ forms a zero-mean, strictly stationary $c$-dependent sequence since for $|s_1 - s_2| > c$, the random variables $Z_{s_1,n}$ and $Z_{s_2,n}$ are based on the two non-overlapping sets of i.i.d. random variables $\{\varepsilon_{t+[u_{s_1}n]-M_n-1}, t = 1, 2, \ldots, m_n\}$ and $\{\varepsilon_{l+[u_{s_2}n]-M_n-1}, l = 1, 2, \ldots, m_n\}$ respectively. Furthermore, for $Var(\widetilde{T}_n)$, we get, for $|s_1 - s_2| \leq c$,

$$Cov(Z_{s_1,n}, Z_{s_2,n}) \to \frac{2}{\pi^2 H_2^4}\left(\int_0^{1-\frac{|s_1-s_2|}{c}} h^2(x)h^2\left(x + \frac{|s_1-s_2|}{c}\right)\mathrm{d}x\right)^2$$

$$\times \int_{-2\pi}^{2\pi}(K*K)(y)\,\mathrm{d}y.$$

Convergence of $N^{-1/2} \sum_{s=1}^{N} Z_{s,n}$ to the desired Gaussian distribution then follows by a CLT for $c$-dependent sequences; see Theorem 6.4.2 of Brockwell and Davis (1991).

The proof of the theorem is completed because straightforward but cumbersome calculations show that $R_{1,n} \to 0$ and $R_{2,n} \to 0$ in probability; see Paparoditis (2006). □

**Proof of Theorem 4.1.** Let $V(u, \lambda) = I_{m_n}(u, \lambda)/g(\lambda) - 1$, $G(u, \lambda) = (f(u, \lambda)/g(\lambda) - 1)$ and write $T_n$ as

$$T_n = \frac{1}{N} \sum_{s=1}^{N} \int (V_n^2(u_s, \lambda) - V^2(u_s, \lambda)) \, d\lambda + \frac{1}{N} \sum_{s=1}^{N} \int (V^2(u_s, \lambda) - G^2(u_s, \lambda)) \, d\lambda$$

$$+ \frac{1}{N} \sum_{s=1}^{N} \int G^2(u_s, \lambda) \, d\lambda.$$

The assertion of the theorem follows because $N^{-1} \sum_{s=1}^{N} \int (V_n^2(u_s, \lambda) - V^2(u_s, \lambda)) \, d\lambda \le O_P(1) \int_{-\pi}^{\pi} |\widehat{g}_h(\lambda) - g(\lambda)| \, d\lambda$, $\int_{-\pi}^{\pi} |\widehat{g}_h(\lambda) - g(\lambda)| \, d\lambda \to 0$ in probability (van der Vaart (1998), Proposition 2.29) and $N^{-1} \sum_{s=1}^{N} \int (V^2(u_s, \lambda) - G^2(u_s, \lambda)) \, d\lambda \le O_P(r_n)$. □

**Proof of Theorem 4.2.** Using the notation $L_n(u, \lambda) = (I_{m_n}(u, \lambda)/f(u, \lambda) - 1)$ and $G_n(u, \lambda) = (f(u, \lambda)/\widehat{g}_h(\lambda) - 1)$ and substituting $I_{m_n}(u, \lambda_j)/\widehat{g}_h(\lambda_j) - 1 = L_n(u, \lambda_j) f(u, \lambda_j)/\widehat{g}_h(\lambda_j) + G_n(u, \lambda_j)$, we get

$$T_n = \frac{1}{m_n^2 N} \sum_{s=1}^{N} \int \left\{ \sum_j K_b(\lambda - \lambda_j) L_n(u_s, \lambda_j) \right\}^2 d\lambda$$

$$+ \frac{1}{m_n^2 N} \sum_{s=1}^{N} \int \left\{ \sum_j K_b(\lambda - \lambda_j) L_n(u_s, \lambda_j) G_n(u_s, \lambda_j) \right\}^2 d\lambda$$

$$+ \frac{1}{m_n^2 N} \sum_{s=1}^{N} \int \left\{ \sum_j K_b(\lambda - \lambda_j) G_n(u_s, \lambda_j) \right\}^2 d\lambda$$

$$+ \frac{2}{m_n^2 N} \sum_{s=1}^{N} \sum_{j_1} \sum_{j_2} \int K_b(\lambda - \lambda_{j_1}) K_b(\lambda - \lambda_{j_2}) \, d\lambda$$

$$\times \{ L_n(u_s, \lambda_{j_1}) L_n(u_s, \lambda_{j_2}) G_n(u_s, \lambda_{j_2}) + L_n(u_s, \lambda_{j_1}) G_n(u_s, \lambda_{j_2})$$

$$+ G_n(u_s, \lambda_{j_1}) L_n(u_s, \lambda_{j_1}) G_n(u_s, \lambda_{j_2}) \}$$

$$= \sum_{i=1}^{6} T_{i,n}$$

with an obvious definition for the $T_{i,n}$'s. To establish the theorem, we show that, in probability,

$$\sqrt{m_n N} T_{1,n} - M_H \frac{\sqrt{N}}{\sqrt{m_n} b} \int_{-\pi}^{\pi} K^2(x)\, dx \to 0, \quad (6.12)$$

$$\sqrt{m_n N} T_{2,n} - M_H \frac{\sqrt{N}}{\sqrt{m_n} b} \int_{-\pi}^{\pi} K^2(x)\, dx \frac{1}{2\pi} \int_0^1 \int_{-\pi}^{\pi} \left( \frac{f(u,\lambda)}{g(\lambda)} - 1 \right)^2 d\lambda\, du \to 0, \quad (6.13)$$

$$\sqrt{m_n N} T_{4,n} \to 0, \quad (6.14)$$

where $M_H = H_4/H_2^2$, and that

$$\sqrt{m_n N}(T_{3,n} + T_{5,n} + T_{6,n} - D_n^2) \Rightarrow \mathcal{N}(0, v^2). \quad (6.15)$$

To see (6.12), note that $\sqrt{m_n N} E(T_{1,n}) = M_H N^{1/2} m_n^{-1/2} b^{-1} \int_{-\pi}^{\pi} K^2(u)\, du + O(N^{1/2} m_n^{-1/2})$. Furthermore, if Assumption 1 is satisfied, then $I_{m_n}(u,\lambda) = \widetilde{I}_{m_n}(u,\lambda) + O(m_n^2 n^{-2})$, where

$$\widetilde{I}_{m_n}(u,\lambda) = (2\pi H_{2,m_n}(0))^{-1} \left| \sum_{t=1}^{m_n} h_{t,m_n} X_{t+[un]-M_n-1,n} \exp\{-\mathrm{i}\lambda t\} \right|^2$$

is the tapered local periodogram of the series $X_1(u), X_2(u), \ldots, X_n(u)$, where $X_t(u) = \sum_j a(u,j)\varepsilon_{t-j}$. By Assumption 6 and simple algebra, we get $\sqrt{m_n N} T_{1,n} = \sqrt{m_n N} \widetilde{T}_{1,n} + \mathrm{o}_P(1)$, where $\widetilde{T}_{1,n}$ is obtained by replacing $I_{m_n}(u,\lambda)$ by $\widetilde{I}_{m_n}(u,\lambda)$ in $T_{1,n}$. Now, $Var(\widetilde{T}_{1,n}) = O(m_n^{-2} N^{-1} b^{-1})$, which implies that $Var(\sqrt{m_n N} \widetilde{T}_{1,n}) = O(m_n^{-1} b^{-1}) \to 0$ as $n \to \infty$.

To establish (6.13), we use $f/\widehat{g}_h - 1 = (f/g - 1) - (g/\widehat{g}_h - 1)f/g$ and write $\sqrt{m_n N} T_{2,n} = \sqrt{m_n N} T_{2,n}^{(1)} + R_{2,n}$, where

$$\sqrt{m_n N} T_{2,n}^{(1)} = \sqrt{\frac{m_n}{N}} \sum_{s=1}^{N} \int \left\{ \frac{1}{m_n} \sum_j K_b(\lambda - \lambda_j) L_n(u_s, \lambda_j) G(u_s, \lambda_j) \right\}^2 d\lambda$$

with $G(u,\lambda) = (f(u,\lambda)/g(\lambda) - 1)$ and the remainder $R_{2,n}$ converges to zero in probability by Lemma 6.1 and Assumption 6 since it is of order $O(m_n^{1/2} N^{1/2} h^4 + N^{1/2}(nbh)^{-1/2} + m_n^{1/2} N^{1/2} b^2 h^2)$. Now, arguing exactly as in the proof of assertion (6.12), it follows that

$$\sqrt{m_n N} T_{2,n}^{(1)} - M_H \frac{\sqrt{N}}{\sqrt{m_n} b} \int_{-\pi}^{\pi} K^2(x)\, dx \frac{1}{2\pi} \int_0^1 \int_{-\pi}^{\pi} \left( \frac{f(u,\lambda)}{g(\lambda)} - 1 \right)^2 d\lambda\, du \to 0$$

and $Var(\sqrt{m_n N} T_{2,n}^{(1)}) = O(m_n^{-1} b^{-1}) \to 0$ as $n \to \infty$.

To show that (6.14) is true, we use the same decomposition for $f/\widehat{g} - 1$ as in (6.13) and verify that $\sqrt{m_n N} T_{4,n} = \sqrt{m_n N} T_{4,n}^{(1)} + R_{4,n}$, where

$$\sqrt{m_n N} T_{4,n}^{(1)} = \frac{2}{\sqrt{m_n^3 N}} \sum_{s=1}^{N} \sum_{j_1} \sum_{j_2} \int K_b(\lambda - \lambda_{j_1}) K_b(\lambda - \lambda_{j_2})$$
$$\times L_n(u_s, \lambda_{j_1}) L_n(u_s, \lambda_{j_2}) G(u_s, \lambda_{j_2}) \, d\lambda,$$

and, for $R_{2,n}$, the remainder $R_{4,n}$ satisfies, by Lemma 6.1 and Assumption 6, $R_{4,n} = O(m_n^{1/2} N^{1/2} h^4 + N^{1/2} (nbh)^{-1/2} + m_n^{1/2} N^{1/2} b^2 h^2) \to 0$ as $n \to \infty$. Now, as for $\sqrt{m_n N} T_{2,n}^{(1)}$, the variance of $\sqrt{m_n N} T_{4,n}^{(1)}$ is of order $O(m_n^{-1} b^{-1})$, while

$$E(\sqrt{m_n N} T_{4,n}^{(1)}) = \frac{\sqrt{N}}{\pi \sqrt{m_n b}} \int K^2(x) \, dx \frac{1}{N} \sum_{s=1}^{N} \int G(u_s, \lambda) \, d\lambda + O(N^{1/2} m_n^{-1/2})$$
$$= O(m_n^{-1/2} N^{-1/2}) + O(N^{1/2} m_n^{-1/2}) \to 0,$$

by Assumption 6, where the last equality follows because $\int_0^1 G(u, \lambda) \, du = 0$.

It remains to establish assertion (6.15). For this, consider the sequence of bivariate random variables $(T_{3,n} - D_n^2, T_{5,n} + T_{6,n})'$. Straightforward calculations yield

$$\sqrt{m_n N}(T_{3,n} - D_n^2, T_{5,n} + T_{6,n})' = \sqrt{m_n N}(J_{1,n}, J_{2,n})' + o_P(1) \qquad (6.16)$$

(see Paparoditis (2006)), where

$$J_{1,n} = -2 \int v(\lambda)(I_n(\lambda) - g(\lambda)) \, d\lambda, \qquad (6.17)$$

$$J_{2,n} = \frac{2}{N} \sum_{s=1}^{N} \int w(u_s, \lambda)(I_{m_n}(u_s, \lambda) - f(u_s, \lambda)) \, d\lambda, \qquad (6.18)$$

$v(\lambda) = \int_0^1 f(u, \lambda)[f(u, \lambda)/g(\lambda) - 1]/g^2(\lambda) \, du$ and $w(u, \lambda) = [f(u, \lambda)/g(\lambda) - 1]/g(\lambda)$.

Using Lemma 6.2, we get $m_n N \, Var(J_{1,n}) \to V_1$ as $n \to \infty$, where

$$V_1 = 16\pi \int v^2(\lambda) g^2(\lambda) \, d\lambda + 8\pi \int \int v(\lambda_1) v(\lambda_2) g_4(\lambda_1, -\lambda_1, \lambda_2) \, d\lambda_1 \, d\lambda_2$$

and $g_4(\lambda_1, -\lambda_1, \lambda_2) = \int_0^1 f_4(u, \lambda_1, -\lambda_1, \lambda_2) \, du$. Furthermore,

$m_n N \, Var(J_{2,n})$

$$= \frac{4m_n}{N} \int \int \sum_{s_1=1}^{N} \sum_{s_2=1}^{N} w(u_{s_1}, \omega) w(u_{s_2}, \lambda)$$
$$\times \{cum(d_n(u_{s_1}, \omega), d_n(u_{s_2}, \lambda)) cum(d_n(u_{s_1}, -\omega), d_n(u_{s_2}, -\lambda))$$

$$+ cum(d_n(u_{s_1}, \omega)d_n(u_{s_2}, -\lambda))cum(d_n(u_{s_1}, -\omega), d_n(u_{s_2}, \lambda))$$

$$+ cum(d_n(u_{s_1}, \omega), d_n(u_{s_1}, -\omega), d_n(u_{s_2}, \lambda), d_n(u_{s_2}, -\lambda))\} \, d\omega \, d\lambda$$

$$= W_{1,n} + W_{2,n} + W_{3,n}$$

and an obvious definition for $W_{i,n}$, $i = 1, 2, 3$. Analyzing each term separately, we get

$$W_{1,n} \to W_1 = \frac{8\pi}{H_2^2} \sum_{|s|<c} \int_{[0,1-\frac{|s|}{c}]} h^2(x)h^2(x+|s|/c) \, dx \int_0^1 \int \widetilde{w}^2(u, \lambda) \, d\lambda \, du,$$

where $\widetilde{w}(u, \lambda) = w(u, \lambda)f(u, \lambda)$. $W_{2,n}$ has the same limit as $W_{1,n}$, while

$$W_{3,n} \to W_3 = \frac{8\pi\kappa_4}{H_2^2} \sum_{|s|<c} \int_{[0,1-\frac{|s|}{c}]} h^2(x)h^2(x+|s|/c) \, dx \int_0^1 \left(\int_{-\pi}^{\pi} w(u, \lambda) \, d\lambda\right)^2 du.$$

Finally, since

$$4m_n \sum_{s=1}^N \int \int v(\lambda_1)w(u_s, \lambda_2) cum(d_n(\lambda_1), d_{m_n}(u_s, \lambda_2)) cum(d_n(-\lambda_1), d_{m_n}(u_s, -\lambda_2)) \, d\lambda_1 \, d\lambda_2$$

$$= \frac{4m_n}{nH_{2,m_n}(0)} \sum_{s=1}^N \int \int v(\lambda_1)w(u_s, \lambda_2)$$

$$\times \sum_{t,l=[u_s n]-M_n+2}^{[u_s n]+M_n+1} h_{t-[u_s n]+M_n-1,m_n} h_{l-[u_s n]+M_n-1,m_n}$$

$$\times A\left(\frac{t}{n}, \lambda_1\right) A\left(\frac{t}{n}, \lambda_2\right) A\left(\frac{l}{n}, -\lambda_1\right) A\left(\frac{l}{n}, -\lambda_2\right)$$

$$\times e^{-i(\lambda_1+\lambda_2)(t-l)} \, d\lambda_d \lambda_2$$

$$+ O(\log(m_n)/m_n)$$

$$\to C_{1,2}^{(1)} = 8\pi \int_0^1 \int v(\lambda)w(u, \lambda)f^2(u, \lambda) \, d\lambda \, du$$

and

$$4m_n \sum_{s=1}^N \int \int v(\lambda_1)w(u_s, \lambda_2) cum(d_n(\lambda_1), d_n(-\lambda_1), d_{m_n}(u_s, \lambda_2), d_{m_n}(u_s, -\lambda_2)) \, d\lambda_1 \, d\lambda_2$$

$$= \frac{4m_n}{nH_{2,m_n}(0)} \sum_{s=1}^N \int \int v(\lambda_1)w(u_s, \lambda_2)$$

$$\times \sum_{t=[u_s n]-M_n+2}^{[u_s n]+M_n+1} h_{t-[u_s n]+M_n-1,m_n}^2 \left| A\left(\frac{t}{n},\lambda_1\right)\right|^2 \left| A\left(\frac{t}{n},\lambda_2\right)\right|^2 + \mathrm{O}(\log(m_n)/m_n)$$

$$\to C_{1,2}^{(2)} = 8\pi \int_0^1 \int \int v(\lambda_1) w(u,\lambda_2) f(u,\lambda_1) f(u,\lambda_1)\, \mathrm{d}\lambda_1\, \mathrm{d}\lambda_2\, \mathrm{d}u,$$

we get $m_n N\, Cov(J_{1,n}, J_{2,n}) \to C_{1,2} = -2C_{1,2}^{(1)} - C_{1,2}^{(2)}$. Using Lemma 6.2 and arguing as in the proof of Theorem 5.10.1 of Brillinger (1981), we get that $\sqrt{m_n N} J_{1,n} \Rightarrow \mathcal{N}(0,V_1)$ as $n \to \infty$. Furthermore, by Theorem A.2 of Dahlhaus (1997), we get $\sqrt{m_n N} J_{2,n} \Rightarrow \mathcal{N}(0,W)$, where $W = 2W_1 + W_3$. An application of the Cramér–Wold device leads to $\sqrt{m_n N}(J_{1,n}, J_{2,n})' \Rightarrow \mathcal{N}((0,0)',V)$, where $V = (v_{r,s})_{r,s=1,2}$ with $v_{1,1} = V_1$, $v_{2,2} = W$ and $v_{1,2} = C_{1,2}$. Approximation (6.16) then yields assertion (6.15) with $v^2 = v_{1,1} + v_{2,2} + 2v_{1,2}$. □

# Acknowledgement

Many thanks are due to a referee for his detailed and helpful comments.